
\ifx\shlhetal\undefinedcontrolsequenc\let\shlhetal\relax\fi
\input amstex
\NoBlackBoxes
\documentstyle {amsppt}
\topmatter
\title {Toward classifying unstable theories \\
Sh500} \endtitle
\author {Saharon Shelah \thanks{\null\newline
Done: \S1: with 457; section 2: 8/92: 2.2 + sufficiency for non
existence of universal; 12/92 rest of \S2 Lecture in Budapest 7/92(?);
\null\newline
Latest revision 25/August/95
\null\newline
Typed at Rutgers: 12/23/94} \endthanks} \endauthor
\affil {Institute of Mathematics \\
The Hebrew University \\
Jerusalem, Israel
\medskip
Rutgers University \\
Department of Mathematics \\
New Brunswick, NJ  USA} \endaffil
\abstract
We prove a consistency results saying, that for a simple (first order)
theory, it is easier to have a universal model in some cardinalities,
than for the theory of linear order. We define additional properties
of first order theories, the $n$-strong order property ($SOP_n$ in short). The
main result is that a first order theory with the 4-strong order property
behaves like linear orders concerning existence of universal models.
\endabstract

\endtopmatter
\document

\def\renewcommand{\newcommand}	       
\edef\cite{\the\catcode`@}%
\catcode`@ = 11
\let\@oldatcatcode = \cite
\chardef\@letter = 11
\chardef\@other = 12
%
%
%
%
\def\@innerdef#1#2{\edef#1{\expandafter\noexpand\csname #2\endcsname}}%
%
%
\@innerdef\@innernewcount{newcount}%
\@innerdef\@innernewdimen{newdimen}%
\@innerdef\@innernewif{newif}%
\@innerdef\@innernewwrite{newwrite}%
%
%
%
\def\@gobble#1{}%
%
%
%
\ifx\inputlineno\@undefined
   \let\@linenumber = \empty 
\else
   \def\@linenumber{\the\inputlineno:\space}%
\fi
%
%
%
\def\@futurenonspacelet#1{\def\cs{#1}%
   \afterassignment\@stepone\let\@nexttoken=
}%
\begingroup 
\def\\{\global\let\@stoken= }%
\\ 
\endgroup
\def\@stepone{\expandafter\futurelet\cs\@steptwo}%
\def\@steptwo{\expandafter\ifx\cs\@stoken\let\@@next=\@stepthree
   \else\let\@@next=\@nexttoken\fi \@@next}%
\def\@stepthree{\afterassignment\@stepone\let\@@next= }%
%
%
%
\def\@getoptionalarg#1{%
   \let\@optionaltemp = #1%
   \let\@optionalnext = \relax
   \@futurenonspacelet\@optionalnext\@bracketcheck
}%
%
%
\def\@bracketcheck{%
   \ifx [\@optionalnext
      \expandafter\@@getoptionalarg
   \else
      \let\@optionalarg = \empty
      \expandafter\@optionaltemp
   \fi
}%
\def\@@getoptionalarg[#1]{%
   \def\@optionalarg{#1}%
   \@optionaltemp
}%
%
%
%
\def\@nnil{\@nil}%
\def\@fornoop#1\@@#2#3{}%
\def\@for#1:=#2\do#3{%
   \edef\@fortmp{#2}%
   \ifx\@fortmp\empty \else
      \expandafter\@forloop#2,\@nil,\@nil\@@#1{#3}%
   \fi
}%
\def\@forloop#1,#2,#3\@@#4#5{\def#4{#1}\ifx #4\@nnil \else
       #5\def#4{#2}\ifx #4\@nnil \else#5\@iforloop #3\@@#4{#5}\fi\fi
}%
\def\@iforloop#1,#2\@@#3#4{\def#3{#1}\ifx #3\@nnil
       \let\@nextwhile=\@fornoop \else
      #4\relax\let\@nextwhile=\@iforloop\fi\@nextwhile#2\@@#3{#4}%
}%
%
%
%
\@innernewif\if@fileexists
\def\@testfileexistence{\@getoptionalarg\@finishtestfileexistence}%
\def\@finishtestfileexistence#1{%
   \begingroup
      \def\extension{#1}%
      \immediate\openin0 =
         \ifx\@optionalarg\empty\jobname\else\@optionalarg\fi
         \ifx\extension\empty \else .#1\fi
         \space
      \ifeof 0
         \global\@fileexistsfalse
      \else
         \global\@fileexiststrue
      \fi
      \immediate\closein0
   \endgroup
}%
%
%
%
%
\def\bibliographystyle#1{%
   \@readauxfile
   \@writeaux{\string\bibstyle{#1}}%
}%
\let\bibstyle = \@gobble
%
%
\let\bblfilebasename = \jobname
\def\bibliography#1{%
   \@readauxfile
   \@writeaux{\string\bibdata{#1}}%
   \@testfileexistence[\bblfilebasename]{bbl}%
   \if@fileexists
      \nobreak
      \@readbblfile
   \fi
}%
\let\bibdata = \@gobble
%
%
\def\nocite#1{%
   \@readauxfile
   \@writeaux{\string\citation{#1}}%
}%
\@innernewif\if@notfirstcitation
%
%
\def\cite{\@getoptionalarg\@cite}%
%
%
\def\@cite#1{%
   \let\@citenotetext = \@optionalarg
   \printcitestart
   \nocite{#1}%
   \@notfirstcitationfalse
   \@for \@citation :=#1\do
   {%
      \expandafter\@onecitation\@citation\@@
   }%
   \ifx\empty\@citenotetext\else
      \printcitenote{\@citenotetext}%
   \fi
   \printcitefinish
}%
\def\@onecitation#1\@@{%
   \if@notfirstcitation
      \printbetweencitations
   \fi
   \expandafter \ifx \csname\@citelabel{#1}\endcsname \relax
      \if@citewarning
         \message{\@linenumber Undefined citation `#1'.}%
      \fi
      \expandafter\gdef\csname\@citelabel{#1}\endcsname{%
\strut
\vadjust{\vskip-\dp\strutbox
\vbox to 0pt{\vss\parindent0cm \leftskip=\hsize 
\advance\leftskip3mm
\advance\hsize 4cm\strut\openup-4pt 
\rightskip 0cm plus 1cm minus 0.5cm ?  #1 ?\strut}}
         {\tt
            \escapechar = -1
            \nobreak\hskip0pt
            \expandafter\string\csname#1\endcsname
            \nobreak\hskip0pt
         }%
      }%
   \fi
   \csname\@citelabel{#1}\endcsname
   \@notfirstcitationtrue
}%
%
%
\def\@citelabel#1{b@#1}%
%
%
\def\@citedef#1#2{\expandafter\gdef\csname\@citelabel{#1}\endcsname{#2}}%
%
%
%
\def\@readbblfile{%
   \ifx\@itemnum\@undefined
      \@innernewcount\@itemnum
   \fi
   \begingroup
      \def\begin##1##2{%
         \setbox0 = \hbox{\biblabelcontents{##2}}%
         \biblabelwidth = \wd0
      }%
      \def\end##1{}
      %
      %
      \@itemnum = 0
      \def\bibitem{\@getoptionalarg\@bibitem}%
      \def\@bibitem{%
         \ifx\@optionalarg\empty
            \expandafter\@numberedbibitem
         \else
            \expandafter\@alphabibitem
         \fi
      }%
      \def\@alphabibitem##1{%
         \expandafter \xdef\csname\@citelabel{##1}\endcsname {\@optionalarg}%
         \ifx\biblabelprecontents\@undefined
            \let\biblabelprecontents = \relax
         \fi
         \ifx\biblabelpostcontents\@undefined
            \let\biblabelpostcontents = \hss
         \fi
         \@finishbibitem{##1}%
      }%
      \def\@numberedbibitem##1{%
         \advance\@itemnum by 1
         \expandafter \xdef\csname\@citelabel{##1}\endcsname{\number\@itemnum}%
         \ifx\biblabelprecontents\@undefined
            \let\biblabelprecontents = \hss
         \fi
         \ifx\biblabelpostcontents\@undefined
            \let\biblabelpostcontents = \relax
         \fi
         \@finishbibitem{##1}%
      }%
      \def\@finishbibitem##1{%
         \biblabelprint{\csname\@citelabel{##1}\endcsname}%
         \@writeaux{\string\@citedef{##1}{\csname\@citelabel{##1}\endcsname}}%
         \ignorespaces
      }%
      %
      %
      \let\em = \bblem
      \let\newblock = \bblnewblock
      \let\sc = \bblsc
      \frenchspacing
      \clubpenalty = 4000 \widowpenalty = 4000
      \tolerance = 10000 \hfuzz = .5pt
      \everypar = {\hangindent = \biblabelwidth
                      \advance\hangindent by \biblabelextraspace}%
      \bblrm
      \parskip = 1.5ex plus .5ex minus .5ex
      \biblabelextraspace = .5em
      \bblhook
      \input \bblfilebasename.bbl
   \endgroup
}%
%
%
\@innernewdimen\biblabelwidth
\@innernewdimen\biblabelextraspace
%
%
%
\def\biblabelprint#1{%
   \noindent
   \hbox to \biblabelwidth{%
      \biblabelprecontents
      \biblabelcontents{#1}%
      \biblabelpostcontents
   }%
   \kern\biblabelextraspace
}%
%
%
%
\def\biblabelcontents#1{{\bblrm [#1]}}%
%
%
\def\bblrm{\rm}%
%
%
\def\bblem{\it}%
%
%
\def\bblsc{\ifx\@scfont\@undefined
              \font\@scfont = cmcsc10
           \fi
           \@scfont
}%
%
%
\def\bblnewblock{\hskip .11em plus .33em minus .07em }%
%
%
\let\bblhook = \empty
%
%
%
\def\printcitestart{[}
\def\printcitefinish{]}
\def\printbetweencitations{, }
\def\printcitenote#1{, #1}
%
%
%
\let\citation = \@gobble
%
%
%
\@innernewcount\@numparams
%
%
\def\newcommand#1{%
   \def\@commandname{#1}%
   \@getoptionalarg\@continuenewcommand
}%
%
%
\def\@continuenewcommand{%
   \@numparams = \ifx\@optionalarg\empty 0\else\@optionalarg \fi \relax
   \@newcommand
}%
%
%
\def\@newcommand#1{%
   \def\@startdef{\expandafter\edef\@commandname}%
   \ifnum\@numparams=0
      \let\@paramdef = \empty
   \else
      \ifnum\@numparams>9
         \errmessage{\the\@numparams\space is too many parameters}%
      \else
         \ifnum\@numparams<0
            \errmessage{\the\@numparams\space is too few parameters}%
         \else
            \edef\@paramdef{%
               \ifcase\@numparams
                  \empty  No arguments.
               \or ####1%
               \or ####1####2%
               \or ####1####2####3%
               \or ####1####2####3####4%
               \or ####1####2####3####4####5%
               \or ####1####2####3####4####5####6%
               \or ####1####2####3####4####5####6####7%
               \or ####1####2####3####4####5####6####7####8%
               \or ####1####2####3####4####5####6####7####8####9%
               \fi
            }%
         \fi
      \fi
   \fi
   \expandafter\@startdef\@paramdef{#1}%
}%
%
%
%
%
\def\@readauxfile{%
   \if@auxfiledone \else 
      \global\@auxfiledonetrue
      \@testfileexistence{aux}%
      \if@fileexists
         \begingroup
            \endlinechar = -1
            \catcode`@ = 11
            \input \jobname.aux
         \endgroup
      \else
         \message{\@undefinedmessage}%
         \global\@citewarningfalse
      \fi
      \immediate\openout\@auxfile = \jobname.aux
   \fi
}%
%
%
\newif\if@auxfiledone
\ifx\noauxfile\@undefined \else \@auxfiledonetrue\fi
%
%
%
%
\@innernewwrite\@auxfile
\def\@writeaux#1{\ifx\noauxfile\@undefined \write\@auxfile{#1}\fi}%
%
%
%
\ifx\@undefinedmessage\@undefined
   \def\@undefinedmessage{No .aux file; I won't give you warnings about
                          undefined citations.}%
\fi
%
%
\@innernewif\if@citewarning
\ifx\noauxfile\@undefined \@citewarningtrue\fi
%
%
%
\catcode`@ = \@oldatcatcode

\newpage

\head {\S0 Introduction} \endhead
\bigskip

Having finished \cite{Sh:a}, an important direction seems to me to try to 
classify unstable theories; i.e. to find meaningful dividing lines.  In
\cite{Sh:10} two such were the strict order property and the independence
property, their disjunction is equivalent to unstability (see
\cite{Sh:a},4.7 = \cite{Sh:c},4.7,p.70).  For theories
without the independence property, we know $S(A)$ (and $S_\triangle(A)$) are
relatively small (see \cite{Sh:10}, Keisler \cite{Ke76}, \cite{Sh:a},
III,\S7,7.3,7.4,II,\S4,4.9,4.10. 
Also for $\lambda > |T|,\{ p \in S(A):p$ does not split over some
$B \subseteq A$ of cardinality $< \lambda\}$ is $\lambda$-dense (see
\cite{Sh:a},7.5 = \cite{Sh:c},7.5,p.140).

Later this becomes interesting in the content of analyzing monadic logic (see
Baldwin Shelah \cite{BlSh:156} representation Baldwin \cite{Bl}.  
By \cite{Sh:197} if ``no monadic expansion of $T$ has the
independence property" is a significant dividing line.

Lately, some model theorist have become interested in finitary versions
called VC dimensions, see Laskowski \cite{Lw92}, Macintyre \cite{Mcxx} (good 
bound for the case of expansion the real field).

More relevant to the present work is the tree property, which is weaker than
the strict order property (in \cite{Sh:c},III,p.171).

In \cite{Sh:93} we try to investigate theories without the tree property, so
called simple.  This can be looked at as a weakening of stable, so: simple
$\Leftrightarrow \kappa_{\text{cdt}}(T) < \infty \Leftrightarrow$ failure
of the tree property $\Leftrightarrow$ suitable local ranks $< \infty$ are
parallel of stable.  We try to do the parallel to (parts of) Ch. II,III of
\cite{Sh:a}, forking being generalized in some ways.  But here instead
showing the number of untrafilters of the Boolean Algebras of formulas
$\varphi(x,\bar a)$ over $A$ is small $(\le |A|^{|T|})$ we show that it can
be decomposed to few subalgebras satisfying a strong chain condition.  In this
context we also succeed to get averages; but the Boolean algebras we get
were derived from normal ones with a little twist.  We did not start with
generalizing the rest of \cite{Sh:a} like supersimple (i.e.
$\kappa_{\text{cdt}}(T) = \aleph_0$, equivalently suitable rank is
$< \infty$).  The test problem in \cite{Sh:93} was trying to characterize the
class of pairs

$$
\align
SP(T) = \biggl\{
(\lambda,\kappa):&\text{ every model of } T \text{ of cardinality }
\lambda \text{ has a} \\
  &\,\kappa \text{-saturated elementary extensions of cardinality } \lambda
\biggr\}.
\endalign
$$
\medskip

\noindent
For simplicity we consider there only 
$\lambda = \lambda^{|T|} > 2^{|T| + \kappa},\kappa > |T|$ and
$(\exists \mu)(\mu = \mu^{< \kappa} \le \lambda \le 2^\mu)$ (if this fails,
see \cite{Sh:576}).  So by \cite{Sh:93} for non-simple $T$, such $(\lambda,
\kappa)$ is
in $SP(T)$ iff $\lambda = \lambda^{< \kappa}$.  If $\mu = \mu^{< \mu} <
\lambda = \lambda^{< \lambda}$, after sutiable forcing preserving $\mu = 
\mu^{< \mu}$, not collapsing cardinals and making $2^\mu = \lambda$, we have
under suitable generalization of MA, so $\kappa < \mu <
\lambda < 2^\mu \Rightarrow (\lambda,\kappa) \in SP(T)$.

It seems much better to use just the cardinal arithmetic assumptions (not
the generalizations of MA).  This call to investigate problem of
${\Cal P}^-(n)$-amalgamation (see \cite{Sh:87b}, \cite{Sh:c},XII,\S5).
For the case of $n = 3$ this means that
\medskip
\roster
\item "{$(*)_3$}"  if $p_0(\bar x,\bar y),p_1(\bar x,\bar z),p_2(\bar y,
\bar z)$, complete types over $A$, each saying the two sequences of variables
are ``independent" in suitable ways (like nonforking) then we extend the
union of the three (preserving ``independence").
\endroster
\medskip

\noindent
Now $(*)_3$ can be proved. \cite{Sh:93},Claim 7.8,p.201,(3.5,p.187).  
But the proof does not work for higher $n$, naturally counterexamples for
the amalgamation should give counterexample to membership in SP.  This was
carried out by finding counterexamples in a \underbar{wider framework}: 
saturation inside $P$ in \cite{Sh:126}; but we could still hope that for 
the ``true" one there is a positive one.

For long, I was occupied elsewhere and not look into it, but eventually
Hrushovski becomes interested (and through him, others) and we try 
to explain the relevnt research below.  Also, it could be asked if simple
unstable theories ``occurs in nature", ``are important to algebraic
applications".  The works cited below forms a positive answer (note 
that, quite natural, those examples concentrate on the lower part of 
hierarchy, like strongly minimal or finite Morley rank).

On the one hand, Hrushovski, continuing \cite{Sh:126}, prove that there
are simple theories with bad behaviour for ${\Cal P}(n)$ so in the result
above the cardinal arithmetic are not enough.

On the other hand, by Hrushovski Pillay \cite{HrPi} in specific cases 
(finite ranks) relevant cases of
$(*)_n$ are proved, for $n > 3$ under very specific conditions: for
$n = 3$ more general; but the relationship with \cite{Sh:93},7.8 of
$(*)$ was not clarified (in both cases the original rank does not work;
the solution in \cite{Sh:93} is to use dnwd (= ``do not weakly divide"),
Hrushovski changes the rank replacing ``contradictory" by having small
rank; this seems a reasonable approach only for supersimple theories and
was carried only for ones with finite rank, and it gives more information
in other respects.

In Hrushovski \cite{Hr1} let 
${\frak C}_0$ be the monster model for a strongly minimal
theory with elimination of imaginaries, $A \subseteq {\frak C}, A =
\text{dcl }A$, such that every $p \in S^m(A)$ with multiplicity $1$ is
finitely satisfiable in $A$, now $\text{Th}({\frak C},A)$ is simple (of
rank $1$) and we can understand PAC in general content.  Hrushovski
\cite{Hr2} does parallel thing for finite rank.
\bigskip

\centerline{$* \qquad * \qquad *$}
\bigskip

We turn to the present work.  First section deals with the existence of
universal models.  Note that existence of saturated models can be 
characterized nicely by stability (see \cite{Sh:c}).

By Kojman Shelah \cite{KjSh:409}, the theory of linear order and more
generally theories without the strict order property has universal models
in ``few" cardinals.

By \cite{Sh:457} give a sufficient conditions for a consistency of ``there
is in $\mu^{++}$ a model of $T$ universal for models of $T$ of cardinality
$\mu^+$", which we use below.

The main aim is to show that all simple theories behave ``better" in this
respect than the theory of linear order.  Specifically, it is consistent
that $\aleph_0 < \lambda = \lambda^{< \lambda},2^\lambda > \lambda^{++}$,
moreover, there is a club guessing $\langle C_\delta:\delta < \lambda^+,
\text{cf}(\delta) = \lambda \rangle$, and every simple $T$ of cardinality
$< \lambda$ has a model in $\lambda^{++}$ universal for $\lambda^+$.  For
this we represent results of \cite{Sh:93} and do the things needed
specifically for the use of \cite{Sh:457}. \newline
See 1.4A(2).
\newpage

\head {\S1 Simple theories have more universal models} \endhead
\bigskip

\noindent
We quote \cite{Sh:457},5.1.
\proclaim{1.1 Lemma}  Suppose
\medskip
\roster
\item "{(A)}"  $T$ is first order, complete for simplicity with elimination
of quantifiers (or just inductive theory with the amalgamation and disjoint
embedding property).
\item "{(B)}"  $K_{ap}$ is a simple $\lambda$-approximation system such
that every $M \in K_{ap}$ is a model of $T$ hence every $M_\Gamma$, where
for $\Gamma \in K_{\text{md}}$ we let \newline
$M_\Gamma = \cup\{M:M \in \Gamma \}$.
\item "{(C)}"  Every model $M$ of $T$ of cardinality $\lambda^+$ can be
embedded into $M_\Gamma$ for some $\Gamma \in K_{\text{md}}$.
\endroster
\medskip

\noindent
Then
\medskip
\roster
\item "{(a)}"  in 4.9 in $V^P$, there is a model of $T$ of cardinality
$\lambda^{++}$ universal for models of $T$ of cardinality $\lambda^+$.
\item "{(b)}"  So in $V^P$, univ$(\lambda^+,T) \le \lambda^{++} <
2^\lambda$.
\endroster
\endproclaim
\bigskip

\demo{Proof}  Straightforward.
\enddemo
\bigskip

\demo{1.1 Fact}  1) Assume $M \prec N,\bar a \in {}^{\omega >}N$, and
$\triangle$ a finite set of formulas possibly with parameters from $M$.
Then there are a formula $\psi(\bar x,\bar b) \in \text{tp}(\bar a,N)$
such that:
\medskip
\roster
\item "{$(*)$}"  for any $\bar a' \in M$ realizing $\psi(\bar x,\bar b)$,
we can find a $\triangle-2$-indiscernible sequence $\langle \bar a_i:i \le
\omega \rangle$ such that: $\bar a_0 = \bar a',\bar a_{\omega + 1} =
\bar a$; hence we can find an indiscernible sequence $\langle a'_2:i <
\omega \rangle$ (in ${\frak C}$) such that the $\triangle$-type of
$\bar a'_0 \char 94 \bar a'_1$ is the same as that of $\bar a' \char 94
\bar a$.
\endroster
\medskip

\noindent
2)  Assume $2^{\theta + |T|} \le \kappa$ and $M \prec N$, moreover
\medskip
\roster
\item "{$\bigotimes$}"  if $A \subseteq M,|A| \le \kappa,\bar a \in
{}^\theta N$ thens ome $\bar a' \in M$ realizes tp$(\bar a,A,N)$.
\endroster
\medskip

\noindent
\underbar{Then} for any $\bar a \in {}^\theta N$ and $B \subseteq M,|B| \le
\theta$, there is $A \subseteq M < |A| \le \kappa$, such that for every
$\bar a' \in {}^\theta M$ realizing tp$(\bar a,A,N)$ there is a sequence
$\langle \bar a_i:i \le \kappa \rangle$ which is 2-indiscernible over
$B,\bar a_0 = \bar a',\bar a_\kappa = \bar a$, hence there is an indiscernible
sequence $\langle \bar a'_i:i < \omega \rangle$ such that $\bar a'_0 \char 94
\bar a'_1$ realizes the same type as $\bar a' \char 94 \bar a$ over $B$.
\enddemo
\bigskip

\demo{Proof}  Obvious [notes on combination set theory]. \newline
1)  Let $\langle p_i:i < k \rangle$ list the possible $\triangle$-types of
sequences of length $\ell g(\bar a) + \ell(\bar a)$, so $k < \omega$.  For
each $p_i,\psi_i(\bar x,\bar b_i) \in \text{tp}(\bar a,M,N)$ such that, if
possible for no $\bar a' \in {}^{\ell g(\bar a)}M$ realizing $\psi_i(\bar x,
\bar b_i)$ do we have: $\bar a' \char 94 \bar a$ realizes $p_i$; (if there..)
\medskip

Now $\psi(\bar x,\bar b) =: \dsize \bigwedge_{i < k} \psi_i(\bar x,\bar b)$
is as required. \newline
2) Similar. \newline
3) $\langle p_i:i < k \rangle$ list the complete $2 \ell g(\bar a)$-types over
$B$.  Use $p_i \subseteq \text{ tp}(\bar a,M,N),|p_i| \le \theta + |T|$
instead $\psi_i(\bar x,\bar b_i),A = \dsize \bigcup_{i < k} \text{ Dom }p_i$.
\enddemo
\bigskip

\proclaim{1.2 Theorem}  If $T$ is a complete simple (f.o.) theory, $|T| <
\lambda$ \underbar{then} $T$ satisfies the assumption of 1.1 (hence its
conclusions).
\endproclaim
\bigskip

\remark{1.2A Remark}  1) We can get results for a theory $T$ of cardinality
$\le \lambda$ under stronger assumptions on $T$. \newline
2) Though not always necessary, in this section we'll assume $T$ is simple.
\newline
3) Also this section is not written in a way focused on Theorem 1.2, but
leisurely relook at simple theories.
\endremark
\bigskip

\demo{Proof}  Without loss of generality $T$ has elimination of quantifiers.

We first represent (in 1.3 - 1.10) the needed definitions and facts on simple
theories from \cite{Sh:93} (adding notation and some facts), then say a little
more and prove the theorem.  So for a while we work in a fixed 
$\bar \kappa$-saturated model ${\frak C}$ of $T, \bar \kappa$ big enough.
So $M,N$ denotes elementary submodels of ${\frak C}$ of cardinality
$< \bar \kappa$, $A,B,C,D$ denote subsets of ${\frak C}$ of cardinality
$< \bar \kappa$ and $\bar a, \bar b, \bar c, \bar d$ denote sequences of
elements of ${\frak C}$ of length $< \bar \kappa$, usually finite.  Let
$\bar a/B = \text{tp}(\bar a,B) = \{ \varphi(\bar x,\bar b):\bar b \in
{}^{\omega >}B,\varphi \in L(T)$ is first order and ${\frak C} \models
\varphi[\bar a,\bar b]\}$.
\enddemo
\bigskip

\definition{1.3 Definition}  1) We say that ``$p(\bar x)$ does not weakly
divide over $(r,B)$" (in short $p$ dnwd (that is does not weakly divide)
over $(r,B)$;
we write over $B$ when $r = p \restriction B$, we write over $r$ if
$B = \text{ Dom}(r)$), where $r = r(\bar x)$ is a type over $B$ (and
$\bar x$ may be infinite) when: if $\bar b \in B$ and $\psi = \psi(\bar x^1,
\dotsc,\bar x^n,\bar y)$ a formula (where $\ell g(\bar x^\ell) = \ell g
(\bar x),\bar x^\ell$ with no repetition, $\langle \bar x^\ell:\ell = 1,n
\rangle \char 94 \langle y \rangle$ pairwise disjoint) and (see definition
1.3(2) below) $[r]^\psi$ is finitely satisfiable (in ${\frak C}$) then so
is $[r \cup p]^\psi$ (see Definition 1.3(2) below). \newline
2)  If $\psi = \psi(\bar x^1,\dotsc,\bar x^n)$ (possibly with parameters),
$q = q(\bar x)$ then \newline
$[q]^\psi = \{ \psi \} \cup \dsize \bigcup^n_{\ell = 1}
q(\bar x^\ell)$. \newline
3) $p(\bar x)$ divides over $A$ if for some formula $\psi(\bar x,\bar a)$
we have $p \vdash q(\bar x,\bar a)$ and for some indiscernible sequence
$\langle \bar a_\ell:\ell < \omega \rangle$ over $A,\bar a = \bar a_0$, and
$\{ \varphi(\bar x,\bar a_i):i < \omega \}$ is $(< \omega)$-contradictory 
where a set $p$ of formulas is $n$-contradictory if for any distinct 
$\varphi_1,\dotsc,
\varphi_n \in p,\{ \varphi_1,\dotsc,\varphi_n\}$ is not realized (in
${\frak C}$), and $(< \omega)$-contradictory means: $n$-contradictory for 
some $n$).  We write dnd for ``does not divide". \newline
4)  The type $p$ fork over the set $A$ if for some $n < \omega$ and formulas
$\varphi_\ell(\bar x,\bar a_\ell)$ for $\ell < n$ we have: $p \vdash \dsize
\bigvee_{\ell < n} \varphi_\ell(\bar x,\bar a_\ell)$ and for each $\ell < n$
the formula $\varphi_\ell(\bar x,\bar a_\ell)$ divides over $A$.  We use
``dnf" as shortening for ``does not fork". \newline
5)  The type $p$ is finitely satisfiable (finitely satisfiable) in $A$ (or
in $\bold I$) if every finite subset $p'$ of $p$ is realized by some sequence
from $A$ (or a member of $\bold I$). \newline
6)  If $D$ is an ultrafilter on Dom$(D) = \bold I$ (where all members of
$\bold I$ have the same length, say $m$) then $Av(B,D) =: \{ \varphi(\bar x,
\bar a):\{\bar b \in \text{ Dom}(D):\varphi(\bar b,\bar a)\} \in D\}$.
\enddefinition
\bigskip

\definition{1.4 Definition}  We say ``$\bar a/A$ (or $\text{tp}(\bar a,A)$)
weakly divides over $B$" if $B \subseteq A$ and $\text{tp}(\bar a,A)$ weakly
divides over $(\text{tp}(\bar a,B),B)$, (similarly for does not weakly
divide).
\enddefinition
\bigskip

\remark{1.4A Remark}  1) An equivalent formulation is ``$\bar a/A$ is an
extension of $\bar a/B$ with the same degree for most $(\triangle,\aleph_0,
k)$"; see 1.5(8) below. \newline
2) On ``divides", ``fork", ``weakly divide" see \cite{Sh:93},Def.1.1,1.2,
2.7(2) respectively.  On the first two also \cite{Sh:a}, but there the focus
is on stable theories.  On ``finitely satisfiable" see \cite{Sh:a},Ch.VII,\S4.
We present here most of their properties, ignoring mainly the connections
with suitable degrees and indiscernibility and the derived Boolean algebras
of formulas (satisfying chain conditions).  For stable $T$ the notions of
Def. 1.3 collapse becoming equivalent (finitely satisfiable - only when the
set is a model, see \cite{Sh:a},Ch.III). \newline
Basic properties are (mostly check directly, but 0A,6,8,9 are quoted).
\endremark
\bigskip

\proclaim{1.5 Claim}  0) [Implications]  If $p$ divides over $A$ then $p$
forks over $A$. \newline
0A)  If $p$ forks over $A$ then $p$ weakly divides over $(p \restriction A,
\emptyset)$, (by \cite{Sh:93},2.11(1),p.184, in its proof we have rely not
only of \cite{Sh:93},2.10(2) + 2.9(3) but also on \cite{Sh:93},2.4(3)).
\newline
0B)  If a type $p$ is finitely satisfiable in $A$ \underbar{then} $p$ does
not fork over $A$. \newline
1)  [Monotonicity]  If $B \subseteq A_1 \subseteq A_2 (\subseteq {\frak C}),
\bar a/A_2$ does not widely divide over $B$. \newline
1A)  If $p$ does not divide over $A,A \subseteq A_1$ and $p_1 \subseteq p$
or at least $p \vdash p_1$ \underbar{then} $p_1$ does not divide over $A_1$.
\newline
1B)  If $p$ does not fork over $A, A \subseteq A_1$ and $p_1 \subseteq p$
or at least $p \vdash p_1$ \underbar{then} $p_1$ does not fork over $A_1$.
\newline
1C)  If $p$ does not widely divide over $(r,A), A \subseteq A_1, r_1 \vdash
r$ and $p_1 \subseteq p$ or at least $p \vdash p_1$ \underbar{then} $p_1$
does not widely divide over $(r_1,A_1)$. \newline
2) [Local character]  $\bar a/A$ does not widely divide 
over $B$ \underbar{iff} for every finite subsequence $\bar a'$ of $\bar a$
and finite subset $A'$ of $A$, $\bar a'/(A' \cup B)$ does not widely divide
over $B$. \newline
2A)  The type $p$ does not weakly divide over $(r,B)$ \underbar{iff} 
every finite $p' \subseteq p$ dnwd over $(r,B)$. \newline
2B)  The type $p$ does not divide over $A$ \underbar{iff} for every finite
$p' \subseteq p$ does not divide over $A$, \underbar{iff} some finite
conjunction $\varphi$ of members of $p$ satisfies the requirement in
Definition 1.3(1). \newline
2C)  The type $p$ does not fork over $A$ \underbar{iff} every finite
$p' \subseteq p$ does not fork over $A$. \newline
3) [More monotonicity]  Assume Rang$(\bar a') = \text{ Rang}(\bar a'')$,
\underbar{then}: $\bar a'/A$ dnwd over $B$ iff
$\bar a''/A$ dnwd over $B$. \newline
3A)  If $B \subseteq A$, Rang $\bar a'' \subseteq \text{ acl}(B \cup
\bar a')$ and $\bar a' \backslash \bar A$ dnwd over $B$
\underbar{then} $\bar a''/A$ dnwd over $B$. \newline
3B) Similarly to 3), 3A) for ``does not divide" and for ``does not fork"
and for ``dnwd over $(r,B)$". \newline
4) [Transitivity] If $A_0 \subseteq A_1 \subseteq A_2$ and $\bar a/
A_{\ell + 1}$ dnwd over $A_\ell$ for $\ell = 0,1$
\underbar{then} $\bar a/A_2$ dnwd over $A_0$. \newline
5) [Extendability] If $B \subseteq A \subseteq A^+,p$ an $m$-type over
$A$ and $p$ does not fork over $B$ \underbar{then} $p$ has an extension
$q \in S^m(A^+)$ which does not fork over $B$ (clear or see 
\cite{Sh:93},2.11(3)). \newline
5A) If $p$ is finitely satisfiable in $A$ and (Dom $p) \cup A \subseteq B$
\underbar{then} we can extend $p$ to a complete type over $B$ finitely
satisfiable in $A$. \newline
6) [Trivial nice behaviour]  $\bar a/A$ does not fork over $A$ (by
\cite{Sh:93},2.11(2)). \newline
6A)  For a set $A$ and an $m$-type $p$ we have: $p$ does not widely divide
over $(p,A)$ (check). \newline
6B) Every $m$-type over $M$ is finitely satisfiable in $M$. \newline
7) [Continuity] If $p_i$ does not widely divide over $(r_i,B_i)$ for
$i < \delta$ and \newline
$i < j < \delta \Rightarrow p_i \subseteq p_j \and r_i \subseteq r_j \and
B_i \subseteq B_j]$ \underbar{then} $\dsize \bigcup_{i < \delta} p_i$
does not widely divide over $(\dsize \bigcup_{i < \delta} r_i, \dsize
\bigcup_{i < \delta} B_i)$. \newline
7A) If $\langle A_i:i < \delta \rangle$ is increasing, $\langle B_i:i <
\delta \rangle$ is increasing and $C/B_i$ is finite satisfiable in $A_i$
for each $i < \delta$ then $C/ \dsize \bigcup_{i < \delta} B_i$ is finite
satisfiable in $\dsize \bigcup_{i < \delta} A_i$. \newline
[Why? if $p \subseteq C/ \dsize \bigcup_{i < \delta} B_i$ is finite then for
some $j$ it is over $B_j$ hence $\subseteq C/B_j$ is satisfiable in $A_j$
hence is satisfiable in $\dsize \bigcup_{i < \delta} A_i$]. \newline
8) [Degree] Let $\bar x_m = \langle x_\ell:\ell < m \rangle,E_m$ be an
ultrafilter on $\Omega_m =: \{(\triangle,k):\triangle = \triangle(\bar x_m)
\subseteq L(T)$ finite, $k < \omega\}$ such that for every $(\triangle_0,
k_0) \in \Omega_m$ the following set belongs to $E_m$:

$$
\{(\triangle,k) \in \Omega_m:\triangle \subseteq \triangle_0 \text{ and }
k_0 < k\}.
$$
\medskip

\noindent
If $p(\bar x)$ is a type over $A,\ell g(\bar x) = m$, \underbar{then} for
some complete type $q(\bar x)$ over $A$ extending $p$ for the $E_m$-majority
of $(\triangle(\bar x),k)$ we have $D(q(\bar x),\triangle,\aleph_0,k) =
D(p(\bar x),\triangle,\aleph_0,k)$ (by \cite{Sh:93},2.2(5),p.182; of course,
we can use infinite $\bar x$). \newline
In such a case we say: $q(\bar x)$ is an $E_m$-nonforking extension of
$p(\bar x)$ or $q(\bar x) \, E_m$-does not fork over $p(\bar x)$.  If
$p(\bar x) \in S^m(A)$ (so $p = q \restriction A)$ we may replace ``over
$p(\bar x)$" by ``over $A$". \newline
9) [Additivity]  If for every $\alpha < \alpha^*$ the type tp$(\bar b^\alpha,
\bar a \cup A \cup \dsize \bigcup_{\beta < \alpha} \bar b^\beta)$ does not
divide over $A \cup \dsize \bigcup_{\beta < \alpha} \bar b^\beta)$
\underbar{then} tp$_*(\dsize \bigcup_{\beta < \alpha^*} \bar b^\beta,A \cup
\bar a)$ does not divide over $A$ (by \cite{Sh:93},1.5,p.181). \newline
10)  [Finitely satisfiable is average]  Let $\ell g(\bar x) = m$ and
$p = p(\bar x)$ a type.  \underbar{Then} $p$ is finitely satisfiable 
in $\bold I$ iff for some ultrafilter $D$ over $\bold I$ we have $p \subseteq
Av(D,\text{Dom }p)$. \newline
11)  If $D$ is an ultrafilter on $\bold I$, \underbar{then} $\text{Av}(D,A)$
belongs to $S^m(A)$ and is a finitely satisfiable in $\bold I$.
\endproclaim
\bigskip

\proclaim{1.6 Claim}  1) [small basis]  If $p \in S^\varepsilon(A)$ and
$B_0 \subseteq A$ \underbar{then} for some $B$ we have:
\medskip
\roster
\item "{$(\alpha)$}"  $B_0 \subseteq B \subseteq A$
\item "{$(\beta)$}"  $|B| \le |\varepsilon| + |T| + |B_0|$
\item "{$(\gamma)$}"  $p$ does not wind downward over $(p \restriction B,B)$.
\endroster
\medskip

\noindent
2)  If $\bar a / (A \cup \bar b)$ dnwd over $A$ and $A \subseteq A^+ \subseteq
ac \ell[A \cup \{ \bar a' \in {\frak C}:\bar a' \text{ realizes } \bar a /
A\}]$ \underbar{then} there is $\bar b'$ (of the same length as $\bar b$) such
that:
\medskip
\roster
\item "{(a)}"  \underbar{if} 
$\bar a' \subseteq A^+$ and $\bar a'/A = \bar a'/A$ \underbar{then}
$\bar b' \char 94 \bar a'/A = \bar b \char 94 \bar a/A$.
\endroster
\medskip

\noindent
3) [weak symmetry]  If $\bar a/(A \cup \bar b)$ dnwd over $A$ and then
$\bar b/(A \cup \bar a)$ dnd over $A$. \newline
4) Assume $A \subseteq B \cap C$ (all $\subseteq {\frak C}$) and $C/B$ is
finitely satisfiable in $A$ hence $A = B \cap C$.  \underbar{Then} $B/C$
dnwd over $(B/A,A)$.
\endproclaim
\bigskip

\demo{Proof}  1) By \cite{Sh:93},3.3,p.186. \newline
2)  By \cite{Sh:93},2.13,p.185 we can get clause (a). \newline
3)  By \cite{Sh:93},2.14,p.185 it dnd. \newline
4)  Straightforward (e.g. use 1.5(10)). \hfill$\square_{1.6}$
\enddemo
\bigskip

\proclaim{1.7 Theorem}  If $M \prec N \prec {\frak C},\|M\| = \mu,\|N\| =
\mu^+,|T| < \kappa,\mu = \mu^{< \kappa}$,\underbar{then} there are
$M^+ \prec N^+$ such that $N \prec N^+,M \prec M^+,\| M^+\| = \mu,\|N^+\|
= \mu^+$ and:
\medskip
\roster
\item "{$(*)_1$}"  if $B \subseteq A \subseteq N,B \subseteq M,|A| < \kappa,
m < \omega,p \in S^m(A)$ and $p$ dnwd over $(p \restriction B,\emptyset)$
\underbar{then} $p$ is realized in $M^+$.
\item "{$(*)_2$}"  if $B \subseteq A \subseteq N,B \subseteq M,|A| < \kappa,
C \subseteq {\frak C},|C| < \kappa$ and $A/(B \cup C)$ dnwd over $B$
\underbar{then} there is $C' \subseteq M^+$ realizing $C/(B \cup A)$.
\endroster
\endproclaim
\bigskip

\demo{Proof}  Clearly we can prove $(*)_1,(*)_2$ separately.  Now $(*)_2$
is immediate from 1.6(2).  As for $(*)_1$, this is proved in \cite{Sh:93},\S4
(read \cite{Sh:93},4.13,4.14,4.15,p.193 there, so we use \cite{Sh:92},Theorem
3.1 which says that a Boolean algebra of cardinality $\lambda^+$ satisfying
the $\kappa$-c.c.,$\lambda^{< \kappa} = \lambda$ is $\lambda$-centered, i.e.
is the union of $\le \lambda$ ultrafilters, so if $\kappa > 2^{|T|}$ we are
done which is enough for our main theorem (1.2 when $\lambda > |T|$).  
Actually repeating the proof of \cite{Sh:93},Theorem 3.1 in the circumstances
of \cite{Sh:93},\S4 show that $\kappa > |T|$ is enough).
\hfill$\square_{1.7}$
\enddemo
\bigskip

\definition{1.8 Definition}  1) $K^0_\lambda$ be 

$$
\align
\biggl\{ 
\bar M:&\bar M = \langle
M_i:i < \lambda^+ \rangle \text{ is } \prec \text{-increasing continuous, 
each } M_i \text{ a model of } T \\
  &\text{ of cardinality } \lambda \text{ and } |M_0| = \emptyset \text{ (we 
stipulate such a model} \\
  &\prec M \text{ for every } M \models T \biggr\}.
\endalign
$$
\medskip

\noindent
2) $\le^0$ is the following partial order on $K^0_\lambda:\bar M^1 \le^0
\bar M^2$ if for $i < j < \lambda^+$ we have $M^1_i \prec M^2_i$.
\enddefinition
\bigskip

\demo{1.9 Observation}  1) If $\langle \bar M^\alpha:\alpha < \delta \rangle$
is an $\le^0$-increasing chain (in $K^0_\lambda$) and $\delta < \lambda^+$
\underbar{then} it has a lub $\bar M:M_i = \dsize \bigcup_{\alpha < \delta}
M^\alpha_i$. \newline
2)  If $M$ is a model of $T$ of cardinality $\lambda^+$, \underbar{then}
for some $\bar M \in K^0_\lambda,M = \dsize \bigcup_{i < \lambda^+} M_i$.
\newline
3)  If $\bar M,\bar N \in K^0_\lambda$ and $\dsize 
\bigcup_{\alpha < \lambda^+} M_\alpha \prec \dsize \bigcup_{\alpha <
\lambda^+} N_\alpha$ \underbar{then} for some club $E$ of $\lambda^+$ for
every $\alpha \in E:M_\alpha \prec N_\alpha$ and $N_\alpha/\dsize
\bigcup_{\beta < \lambda^+} M_\beta$ dnwd over $M_\beta$.
\enddemo
\bigskip

\demo{Proof}  1) Immediate. \newline
2) Use 1.5(2) + 1.6(1). \hfill$\square_{1.9}$
\enddemo
\bigskip

\noindent
Using 1.7 and 1.5(5A) $\lambda^+ \times \lambda$ we get
\demo{1.10 Observation}  Assume $\lambda = \lambda^{< \kappa}$.  For every
$\bar M \in K^0_\lambda$ there is an $\le^0$-increasing continuous sequence
$\langle \bar N_\zeta:\zeta \le \lambda \rangle$, in $K^0_\lambda$,
(so $\bar N_\zeta = \langle N^\zeta_\alpha:\alpha < \lambda^+ \rangle),
\bar N_0 = \bar M$ such that (letting $N_\zeta = \dsize \bigcup_{\alpha
< \lambda^+} N^\zeta_\alpha)$ and [fixing $\chi$, letting
$\bar a_{\zeta,\alpha}$ an enumeration of $|N^\zeta_\alpha|$ of length
$\lambda$] we can add: every type definable in $(H(\chi),\in,<^*_\chi)$ from
$\langle \bar N_\varepsilon:\varepsilon \le \zeta \rangle, \langle
N^{\zeta + 1}_\alpha:\alpha \le \beta \rangle, \langle \bar a_{\varepsilon,
\alpha}:\alpha < \lambda^+ \rangle:\varepsilon \le \zeta \rangle,
\langle \bar a_{\zeta,\alpha}:\alpha \le \beta \rangle$ and
$\langle \alpha_{\zeta,\alpha,j}:j < i \rangle$ and finitely many ordinals
$< \lambda$ is realized in $N^{\zeta + 1}_{\beta + 1}$, hence:
\medskip
\roster
\item "{$(*)_1$}"  if $\alpha < \lambda^+, \zeta \le \lambda,
\text{cf}(\zeta) \in \{ \lambda,1\}, \text{cf}(\alpha) \in \{ \lambda,1\}$,
\newline
$B \subseteq A \subseteq N_{\zeta-1}, B \subseteq N^{\zeta-1}_\alpha,
|A| < \kappa, p \in S^m(A)$ and $p$ dnwd over $(p \restriction B,\emptyset)$
\underbar{then} $p$ is realized in $N^\zeta_\alpha$. \newline
(Note: $\lambda -1 = \lambda)$.
\item "{$(*)_2$}"  if $\alpha \le \beta < \gamma < \lambda^+, \gamma$ is
non-limit, $B \subseteq A \subseteq A^+ \subseteq M_\gamma, |A^+| < \kappa,
\bar a \in M_\gamma,\text{tp}(\bar a,A \cup M_\alpha)$ dnwd over $B$ \underbar
{then} either for some $\bar a' \in N_\gamma$ we have:
$\text{tp}(\bar a,A^+ \cup M_\alpha) = \text{tp}(\bar a',A \cup M_\alpha)$
and tp$(\bar a',A \cup N_\beta)$ dnwd over $B$ or there is no such
$\bar a' \in {\frak C}$.
\item "{$(*)_3$}"  for $\gamma$ non-limit, $\zeta \le \lambda,\text{cf}(\zeta)
\in \{ 1,\lambda\}$ we have: \newline
$N^\zeta_\gamma$ is $\kappa$-saturated (so when $\kappa = \lambda$ it is
saturated).
\endroster
\enddemo
\bigskip

\definition{1.11 Definition}  Let $A,B,C$ be given $(\subseteq {\frak C})$.
\newline
0)  $A \le^{-1}_B C$ means that for every $\bar b \subseteq B, \bar b/(A \cup
C)$ dnwd over $({\frac {\bar b}A},A)$. \newline
1)  $A \le^0_B C$ means that for every $\bar c \subseteq C,
\bar c/(A \cup B)$ dnwd over $({\frac {\bar c}A},\emptyset)$. \newline
2)  $A \le^1_B C$ means there is an increasing continuous sequence
$\langle A-\alpha:\alpha \le \beta \rangle$ such that: $A = A_0, A \cup C =
A_\beta$ and

$$
\gather
\alpha \text{ an even ordinal } < \beta \Rightarrow A_\alpha \le^0_B
A_{\alpha + 1} \\
\alpha \text{ an odd ordinal } < \beta \Rightarrow A_\alpha \le^{-1}_B
A_{\alpha + 1}.
\endgather
$$
\medskip

\noindent
3)  $A \le^2_B C$ means that for some $C', C \subseteq C'$ and $A \le^1_B C'$.
\newline
4)  $A \le^3_B C$ means that for some increasing continuous sequence
$\langle A_\alpha:\alpha \le \beta \rangle$ we have $A = A_0, A \cup C =
A_\beta$ and $A_\alpha \subseteq^2_B A_{\alpha + 1}$.
\enddefinition
\bigskip

\proclaim{1.12 Claim}  0) $A \le^e_B A$ (for $e = -1,0,1,2,3)$. \newline
1)  $A \le^e_B C$ iff $A \le^e_{A \cup B} A \cup C$ (for
$e = -1,0,1,2,3$). \newline
[Why?  For $e = -1$ by 1.5(3A) for $e = 0$ by the definition of $\le^0_B$;
then continued]. \newline
2)  If $A \subseteq B_1 \subseteq B \cup A$ and $A \le^e_B C$ then
$A \le^e_{B_1} C$ (for $e = -1,0,1,2,3,$). \newline
[Why?  By part (1) and: for $e = -1$ trivially, for $e = 0$ by 1.5(1C) and for
$e = 1,2,3$ use earlier cases]. \newline
3) For $e = 1,3$ we have: $\le^e_B$ is a partial order. \newline
[Why?  Read their definition]. \newline
4)  If $e = 1,3$ and $\langle A_\alpha:\alpha \le \beta \rangle$ is
increasing continuous and $A_\alpha \le^e_B A_{\alpha + 1}$ for
$\alpha < \beta$ \underbar{then} $A_0 \le^e_B A_\beta$. \newline
[Why?  Check]. \newline
5)  For $(e^1,e^2) \in \{ (-1,1),(0,1),(1,2),(2,3)\}$, we have:
$A \le^{e^1}_B C$ implies $A \le^{e^2}_B C$. \newline
[Why?  Read]. \newline
6)  If for every $\bar b \subseteq B,\bar b/(A \cup C)$ is finitely
satisfiable in $A$ \underbar{then} $A \le^0_B C$. \newline
[Why?  By 1.6(4) (and Definition 1.11(1))]. \newline
7)  If $A \le^2_B C$ and $C' \subseteq C$ \underbar{then} $A \le^2_B C'$.
\newline
[Why? Read Definition 1.11(2)]. \newline
8)  $A \le^0_B C$ iff $A \le^{-1}_C B$. \newline
[Why?  Read the definitions].
\endproclaim
\bigskip

\proclaim{1.13 Claim}  Let $M \prec N$ and $M \subseteq A$.  Then the
following are equivalent:
\medskip
\roster
\item "{(a)}"  $M \le^3_N A$
\item "{(b)}"  there are $M_0 \prec M_1 \prec M_2$ such that:
{\roster
\itemitem{ (i) }  $M = M_0$
\itemitem{ (ii) }  the type $\text{tp}_*(N,M_1)$ is finitely satisfiable
in $M_0$ and the type \newline
$\text{tp}_*(M_2,M_1 \cup N)$ is finitely satisfiable in $M_1$
\itemitem{ (iii) }  for some elementary map $f,f(A) \subseteq M_2$ and
$f \restriction N$ = identity
\endroster}
\item "{(c)}"  like (b) with $\| M_2 \| \le |T| + |A|$
\item "{(d)}"  $M \le^2_N A$.
\endroster
\endproclaim
\bigskip

\remark{13A Remark}  1) Clause $(ii)$ of $(b)$ implies $M_0 \le^0_n M_1
\le^{-1}_N M_2$. \newline
2)  An equivalent formulation of $(b)$ is
\medskip
\roster
\item "{$(b)'$}"  for some $M_0,M_1,M_2,f$ we have
$M = M_0 \le_{f(N)} M_1 \le_{f(N)} M_2, f \restriction M_0 = \text{id}
_{M_0}, f(A) \subseteq M_2$.
\endroster
\medskip

\noindent
3)  Another formulation is
\medskip
\roster
\item "{$(b)''$}"  like $(b)'$ but $f = \text{id}_{A \cup N}$.
\endroster
\endremark
\bigskip

\demo{Proof}  \underbar{$(c) \Rightarrow (b)$}  Trivial.
\medskip

\noindent
\underbar{$(b) \Rightarrow (c)$}.
\medskip

\noindent
By the Lowenkeim Skolem argument.
\medskip

\noindent
\underbar{$(b) \Rightarrow (d)$}
\medskip

By 1.12(6) clearly $M_0 \le^0_N M_1$ and similarly $M_1 \le^0_{M_2} N$,
hence by 1.2(8) we have $M_1 \le^{-1}_N M_2$.  Hence by 1.12(5),
$M_e \le^1_N M_{e+1}$ (for $e = 0,1$), so by 1.12(3) $M_0 \le^1_N M_2$ hence
by Definition 1.11(3) (and clause $(iii)$ of 1.13(b)), $M = M_0 \le^2_N$
as required.
\medskip

\noindent
\underbar{$(d) \Rightarrow (a)$}  Trivial (by 1.12(5)).
\medskip

So the only (and main) part left is:
\medskip

\noindent
\underbar{$(a) \Rightarrow (b)$}  We know $M \le^3_N A$, by 1.12(1) without
loss of generality $M \subseteq A$, hence there is an increasing continuous
sequence $\langle A_\varepsilon:\varepsilon \le \zeta \rangle$ such that:
$A_0 = M, A_\zeta = A$ and $A_\varepsilon \le^2_N A_{\varepsilon + 1}$.  By
the Definition of $\le^2_n,\le^1_N$ there is an increasing continuous
sequence $\langle B_{\varepsilon + i}:i \le i_\varepsilon \rangle$ such that
$B_{\varepsilon,0} = A_\varepsilon, A_{\varepsilon + 1} \subseteq
B_{\varepsilon,i_\varepsilon}$ and
$B_{\varepsilon,i} \le^{\ell(\varepsilon)}_N B_{\varepsilon,i + 1}$ (where
for $i < i_\varepsilon$ we have $\ell(\varepsilon) \in \{-1,0\}$ and
$\varepsilon = \ell(\varepsilon) \text{ mod } 2$).  Let
$\theta = 2^{|T|} + |N| + \dsize \sum_{\varepsilon < \zeta}
(|i_\varepsilon|) + |B_{\varepsilon,i_\varepsilon}|)^+$ and choose regular
$\mu = \mu^\theta$. 
\medskip

We choose by induction on $\alpha < \mu^+,M_\alpha,N_\alpha$ such that:
medskip
\roster
\item "{(i)}"  $M_\alpha \prec {\frak C}$ is increasing continuous
\item "{(ii)}"  $\| M_\alpha\| = \mu,M \subseteq M_0$,
\item "{(iii)}"  $f_\alpha$ is an elementary mapping,$\text{Dom}(f_\alpha) =
N,\text{Rang}(f_\alpha) = N_\alpha$ and $f_\alpha \restriction M =
\text{ id}_M$
\item "{(iv)}"  $\text{tp}_*(N_\alpha,M_\alpha)$ is finitely satisfiable
in $M$
\item "{(v)}" $N_\alpha \subseteq M_{\alpha + 1}$
\item "{(vi)}"  $M_{\alpha + 1}$ is $\theta^+$-saturated.
\endroster
\medskip

There is no problem to carry the definition.  (First choose $M_\alpha$: if
$\alpha = 0$ to satisfy $(i) + (ii)$, if $\alpha$ is a limit ordinal, as
$\dsize \bigcup_{\beta < \alpha} M_\beta$, and if $\alpha = \beta + 1$ to
satisfy $(i) + (ii) + (v)$.  Second choose $f_\alpha,N_\alpha$ satisfying
$(iii) + (iv)$ which exists by 1.5(10) + (11)).

By using 1.7, $\lambda^+$ times we can find $\bar M^+ = \langle M^+_\alpha:
\alpha < \lambda^+ \rangle$ such that:
\medskip
\roster
\item "{$(A)$}"  $\bar M^+$ is an increasing continuous sequence of
elementary submodels of ${\frak C}$
\item "{$(B)$}"  $\| M^+_\alpha\| \le \mu,M_\alpha \prec M^+_\alpha$
\item "{$(C)_1$}"  if $\alpha < \beta < \mu^+,B_1 \subseteq M_\alpha$ and
$\text{cf}(\alpha) = \mu,B_1 \subseteq B_2 \subseteq M_\beta,|B_2| <
\theta,C \subseteq {\frak C}$ and $C/B_2$ dnwd over $(C/B_2,\emptyset)$
(equivalently, for every finite $\bar c \subseteq C,\bar c/B_2$ dnwd over
$(\bar c/B_1,\emptyset)$) \underbar{then} $C/B_2$ is realized in $M_\alpha$
\item "{$(C)_2$}"  similarly, but we replace the dnwd assumption by
``$B_2/(C_1 \cup C)$ dnwd over $(B_2/B_1,B_1)$". \newline
[Note: we use $(*)_1$ from 1.7 for $(C)_1$ and $(*)_2$ from 1.7 for
$(C)_2$].
\endroster
\medskip

Now let $M = \dsize \bigcup_{\alpha < \mu} M_\alpha,M^+ = \dsize \bigcup
_{\alpha < \mu^+} M^+_\alpha$; and let \newline
$E = \{ \delta < \mu^+:\delta \text{ a limit ordinal and } (M^+_\delta,
M_\delta) \prec (M^+,M)\}$.  Clearly $E$ is a club of $\mu^+$ and
\medskip
\roster
\item "{$(*)$}"  $\delta \in E \Rightarrow \text{ tp}(M^+_\delta,M)$ is
finitely satisfiable in $M_\delta$.
\endroster
\medskip

Choose $\delta \in E$ of cofinality $\mu$.  Now we choose $g_\varepsilon$
by induction on $\varepsilon \le \zeta$ such that:
\medskip
\roster
\item "{$(\alpha)$}"  $g_\varepsilon$ an elementary mapping
\item "{$(\beta)$}"  $\text{Dom}(g_\varepsilon) = N \cup A_\varepsilon$
\item "{$(\gamma)$}"  $g_\varepsilon$ is increasing continuous in
$\varepsilon$
\item "{$(\delta)$}"  $g_\varepsilon \restriction N = f_\delta$
\item "{$(\varepsilon)$}"  $\text{Rang}(g_\varepsilon \restriction
A_\varepsilon) \subseteq M^+_\delta$.
\endroster
\medskip

If we succeed, then we get the desired conclusion (i.e. prove clause $(b)$).
\newline
[Why?  First note that in clause $(b)$ we can omit $f \restriction N =
\text{id}$ by $f \restriction M =$ the identity if in clause (ii) we use
$f(N)$; we call this $(b)'$.  Now $(b)'$ holds with $M,M_\delta,M^+_\delta,
g_\zeta$ here standing to $M_0,M_1,M_2,f$ there).  So it is enough to carry
the induction on $\varepsilon$.  For $\varepsilon = 0$ let $g_\varepsilon =
f_\delta$, and for $\varepsilon$ a limit ordinal let $g_\varepsilon =
\dsize \bigcup_{\xi < \varepsilon} g_\xi$; lastly for $\varepsilon$ a
successor ordinal say $\varepsilon = \xi + 1$, we choose $g_{\varepsilon,i}$
by induction on $i \le i_\varepsilon$ such that:
\medskip
\roster
\item "{$(\alpha)'$}"  $g_{\varepsilon,i}$ an elementary mapping
\item "{$(\beta)'$}"  Dom$(g_{\varepsilon,i}) = N \cup B_{\varepsilon,i}$
\item "{$(\gamma)'$}"  $g_{\varepsilon,i}$ is increasing continuous in $i$
\item "{$(\delta)'$}"  $g_{\varepsilon,0} = g_\varepsilon$
\item "{$(\varepsilon)'$}"  Rang$(g_\varepsilon) \subseteq M^+_\delta$.
\endroster
\medskip

If we succeed then $g_{\varepsilon,i_\varepsilon} \restriction
A_{\varepsilon + 1}$ is as required.  So it is enough to carry the induction
on $i$.  For $i = 0$ let $g_{\varepsilon,i} = g_\varepsilon$, for $i$ limit
let $g_{\varepsilon,i} = \dsize \bigcup_{j < i} g_{\varepsilon,j}$ and for
$i$ a successor ordinal say $j + 1$, use clause $(C)_1$ in the choice of
$M^+_\alpha$ if $j$ even, remembering Definition 1.11(1) and use clause
$(C)_2$ in the choice of $M^+_\alpha$ if $j$ is odd remembering Definition
1.11(0). \hfill$\square_{1.13}$
\enddemo
\bigskip

\proclaim{1.13A Claim}  If $M \prec N, M \le^2_N A_\ell$ for $\ell = 1,2$
then there are $M^+,f_1,f_2$ such that: $M \prec M^+,M \le^2_N M^+$ and for
$\ell = 1,2 f_\ell$ is an elementary mapping, $\text{Dom}(f_\ell) = N \cup
A_\ell,f_\ell \restriction N = \text{ id}_N,f_\ell(A_\ell) \subseteq M^+$.
\endproclaim
\bigskip

\demo{Proof}  Same proof as 1.13 (just shorter).
\enddemo
\bigskip

\definition{1.14 Definition}  1) Let $K^{pr}_0 = \{(M,N):M \prec N <
{\frak C}\}$ and $(M_1,N_1) \le^* (M_2,N_2)$ \underbar{iff} ($(M_e,N_e) \in
K^{pr}_0$ for $e = 1,2$ and $M_1 \prec M_2,N_1 \prec N_2$ and $M_1 \le^2
_{N_1} M_2$ (equivalently, $M_1 \le^3_{N_1} M_2$ (by 1.13)). \newline
2) We define $(M_1,N_1) \le_{fs} (M_2,N_2)$ similarly replacing
``$M_1 \le^2_{N_1} M_2$" by ``$N_1/M_2$ is finitely satisfiable in $M_1$".
\enddefinition
\bigskip

\proclaim{1.15 Claim}  1) $\le^*$ is a partial order on $K^{pr}_0$. \newline
2)  If $\langle(M_\alpha,N_\alpha):\alpha \le \beta \rangle$ is increasing
continuous and $(M_\alpha,N_\alpha) \le^* (M_{\alpha + 1},N_{\alpha + 1})$ for
$\alpha < \beta$ \underbar{then} $(M_0,N_0) \le^* (M_\beta,N_\beta)$. \newline
3)  If $M \prec N$ and $M \le^2_N A$ then for some $(M_1,N_1)$ we have:
$A \subseteq M_2$ and $(M,N) \le^* (M_1,N_1) \in K^{pr}_0$.
\endproclaim
\bigskip

\demo{Proof}  1) If $(M_0,N_0) \le^* (M_1,N_1) \le^* (M_2,N_2)$ then
\medskip
\roster
\item "{(i)}"  $M_0 \subseteq M_1 \subseteq M_2$ and $N_0 \subseteq N_1
\subseteq N_2$
\item "{(ii)}"  $M_0 \le^2_{N_0} M_1$ and
\item "{(iii)}"  $M_1 \le^2_{N_1} M_2$
\item "{(iv)}"  $M_\ell \subseteq N_\ell$.
\endroster
\medskip

\noindent
By 1.12(2) + (1) and clause (iii) above
\medskip
\roster
\item "{(v)}"  $M_1 \le^2_{N_0} M_2$,
\endroster
\medskip

\noindent
by 1.12(5) we have (by (ii) and (v) respectively)
\medskip
\roster
\item "{$(ii)'$}"  $M_0 \le^3_{N_0} M_1$
\item "{$(v)'$}"  $M_1 \le^3_{N_0} M_2$
\endroster
\medskip

\noindent
hence by 1.13
\medskip
\roster
\item "{(vii)}"  $M_0 \le^2_{N_0} M_2$
\endroster
\medskip

\noindent
hence $(M_0,N_0) \le^* (M_2,N_2)$ holds by (i), (iv) and (vii). \newline
2)  Similarly using 1.12(4) + (1.13). \newline
3)  Use 1.13 (see 1.13a(3)) so there are $M_0 \prec M_1 \prec M_2$ such that
$M = M_0,N/M_1$ fs in $M_0,M_2/(M_1 \cup N)$ fs in $M_1$ and $A \subseteq
M_2$.  So by 1.x $M_0 \le^{-1}_N M_1,M_1 \le^0_N M_1$ hence (see 1.12(x)
$M_0 \le^3_N M_0 \le^3_N M_1$ hence (see 1.12(4)) $M_0 \le^3_N M_2$ hence
for any $N^*,M_2 \cup N \subseteq N^* \prec {\frak C}$ we have
$(M,n) \le^* (M_2,N^*) \in K^{pr}_0$. \hfill$\square_{1.15}$
\enddemo
\bigskip

\definition{1.16 Definition}  $K^{pr}_2 = \{(M,N):\text{the pair } (M,N) \in
K^{pr}_0 \text{ and if } (M,N) \le^* (M',N') \in K^{pr}_0$ then
$M'/N$ is fs in $M\}$.  
\enddefinition
\bigskip

\proclaim{1.17 Claim}  If $(M,N) \in K^{pr}_0$ 
then for some $(M',N')$ we have:
\medskip
\roster
\item "{(a)}"  $(M,N) \le (M',N') \in K^{pr}_0$
\item "{(b)}"  $\| N'\| \le \| N \| + |T|$
\item "{(c)}"  $(M',N') \in K^{pr}_1$ i.e. \newline
if $(M',N') \le^* (M'',N'')$ then $M''/N'$ is fs in $M'$.
\endroster
\endproclaim
\bigskip

\demo{Proof}  Let $\mu = \| N \| + |T|$, assume the conclusion fails.  We
now choose by induction on $\alpha < \mu^+,(M_\alpha,N_\alpha)$ such that:
\medskip
\roster
\item "{(i)}"  $(M_0,N_0) = (M,N)$
\item "{(ii)}"  $(M_\alpha,N_\alpha) \in K^{pr}_0,\| N_\alpha\| \le \mu$
\item "{(iii)}"  $\beta < \alpha \Rightarrow (M_\beta,N_\beta) \le^*
(M_\alpha,N_\alpha)$
\item "{(iv)}"  for limit $\delta$ we have $(M_\delta,N_\delta) =
\left( \dsize \bigcup_{\alpha < \delta} M_\alpha,\dsize \bigcup_{\alpha <
\delta} N_\delta \right)$
\item "{(v)}"  $M_{\alpha + 1}/N_\alpha$ is not fs in $M_\alpha$.
\endroster
\medskip

\noindent
For $\alpha = 0$ see (i) for $\alpha$ limit see (iv) and 1.15(2) if
$\alpha = \beta + 1$ find $(M_\alpha,N_\alpha)$ satisfying $(M_\beta,N_\beta)
\le^* (M_\alpha,N_\alpha) \in K^{pr}_1$ and satisfying (v).  By Lowenheim
Skolem argument without loss of generality $\| N_\alpha \| \le \mu$ and by
1.15(1) also clause (iii) holds.  For a club of $\delta < \mu^+$ we get
contradiction to clause (v). \hfill$\square_{1.17}$
\enddemo
\bigskip

\demo{1.18 Fact}  1) If $(M,N) \in K^{pr}_1$ and $(M',N') \in K^{pr}_0$
and $(M,N) \le^* (M',N')$ then $(M,M') \le^* (N,N')$. \newline
2)  If $(M,N) \in K^{pr}_1$ and $M \le^2_N A$ then $A/N$ is fs in $M$.
\enddemo
\bigskip

\demo{Proof}  1) By 1.17 we know 
$M'/N$ is fs in $M$ hence by 1.13, $(b) \Rightarrow (d)$ we
know $M \le^2_{M'} N$ which give the desired conclusion. \newline
2)  By 1.13A.
\enddemo
\bigskip

\proclaim{1.19 Claim}  1) If $(M_\alpha,N_\alpha) \in K^{pr}_1$ for
$\alpha < \delta$ and $\langle(M_\alpha,N_\alpha):\alpha < \delta \rangle$
is \newline
$<^*$-increasing \underbar{then} for $\alpha < \delta$

$$
(M_\alpha,N_\alpha) \le^* \left( \dsize \bigcup_{i < \delta} M_i,
\dsize \bigcup_{i < \delta} N_i \right) \in K^{pr}_1.
$$
\bigskip

\noindent
2)  If $(M_\alpha,N_\alpha) \in K^{pr}_1$ and $\langle(M_\alpha,N_\alpha):
\alpha \le \delta \rangle$ is $<^*$-increasing then \newline
$\left( \dsize \bigcup_{\alpha < \delta} M_\alpha,\dsize \bigcup
_{\alpha < \delta} N_\alpha \right) \in K^{pr}_1$ and
$\left( \dsize \bigcup_{\alpha < \delta} M_\alpha,\dsize \bigcup
_{\alpha < \delta} N_\alpha \right) \le^* (M_\delta,N_\delta)$.
\endproclaim
\bigskip

\demo{Proof}  We prove both together by induction on $\delta$. \newline
0)  By the induction hypothesis without loss of generality $\langle
(M_\alpha,N_\alpha):\alpha < \delta \rangle$ is increasing continuous.
\newline
1)  Clearly $(M_\alpha,N_\alpha) \le^* \left( \dsize \bigcup_{i < \delta} M_i,
\dsize \bigcup_{i < \delta} N_i \right) \in K^{pr}_0$ (see 1.17(2)).
Suppose \newline
$\left( \dsize \bigcup_{i < \delta} M_i,\dsize \bigcup
_{i < \delta} N_i \right) \le^* (M,N)$.  So by 1.17(1), for $\alpha < \delta,
(M_\alpha,N_\alpha) \le^* (M,N)$, but $(M_\alpha,N_\alpha) \in K^{pr}_1$ hence
$M/N_\alpha$ is fs in $M_\alpha$.  But this implies $M / \dsize \bigcup
_{\alpha < \delta} N_\alpha$ is fs in $\dsize \bigcup_{\alpha < \delta}
M_\alpha$ by 1.5(7A). \newline
2)  As we are proving by induction on $\delta$; without loss of generality
$\langle (M_\alpha,N_\alpha):\alpha < \delta \rangle$ is $\le^*$-increasing
continuous, so by part (1), $(M_\alpha,N_\alpha) \le^*
\left( \dsize \bigcup_{i < \delta} M_i,\dsize \bigcup_{i < \delta} N_i
\right) \in K^{pr}_1$ for $\alpha < \delta$.  Now for $\alpha < \delta,
(M_\alpha,N_\alpha) \le^* (M_\delta,N_\delta)$ and $(M_\alpha,N_\alpha) \in
K^{pr}_1$ clearly $M_\delta/N_\alpha$ is fs in $M_\alpha$ hence by
1.5(7A), $M_\delta / \dsize \bigcup_{\alpha < \delta} N_i$ is fs in
$\dsize \bigcup_{i < \delta} M_i$, hence by 1.13 $\dsize \bigcup_{i < \delta}
M_i \le^3_{\dsize \bigcup_{\alpha < \delta} N_i} \dsize \bigcup_{i < \delta}
N_i$ hence $\left( \dsize \bigcup_{i < \delta} M_i,\dsize \bigcup_{i < \delta}
N_i \right) \le^* (M_{\delta + 1},N_{\delta + 1})$.  \hfill$\square_{1.19}$
\enddemo
\bigskip

\noindent
Now we want to apply 1.1.  Toward this (for $\lambda$ as there) we define:
\definition{1.20 Definition}  1) $K^0_{ap} = K^0_{ap} = K^0_{ap}[T] =
K^0_{ap}[T,\lambda]$ is the set of models $M$ of $T$ with universe
$\subseteq \lambda^+$ and cardinality $< \lambda$ such that: $M \cap \lambda
\ne \emptyset$ and $0 < \alpha < \lambda^+$ implies $M \restriction (\lambda
\times \alpha) \prec M$.  For such $M$ let Dom$(M) = \{ \alpha < \lambda^+:
[\lambda \times \alpha, \lambda \times \alpha + \lambda) \cap M \ne \emptyset
\}$.  We now define $<_{K^0_{ap}}$ by: $M \le_{K^0_{ap}} N$ if (both are
in $K^0_{ap}$ and $M \prec N$ and): for every $\alpha \in (0,\lambda^+),
M \restriction (\lambda \times \alpha) <^2_{M \restriction (\lambda \times
\alpha + \lambda)} N \restriction (\lambda \times \alpha)$.
\enddefinition
\bigskip

\demo{1.21 Observation}  So $M \le_{K^0_{ap}} N$ \underbar{iff} both are
in $K^0_{ap},M < N$ and for $\alpha \in (0,\lambda^+)$ we have
$(M \restriction (\lambda \times \alpha),M \restriction (\lambda \times
\alpha + \lambda)) \le^* (N \restriction (\lambda \times \alpha),N 
\restriction (\lambda \times \alpha + \lambda))$.
\enddemo
\bigskip

\proclaim{1.22 Claim}  1) $\le_{K^0_{ap}}$ is a partial order on $K^0_{ap}$.
\newline
[Why?  By 1.15(1)]. \newline
2)  If $\langle M_i:i < \delta \rangle$ is $\le_{K^0_{ap}}$-increasing,
$\dsize \sum_{i < \delta} \| M_i \| < \lambda$ \underbar{then} $M_i
\le_{K^0_{ap}} \dsize \bigcup_{j < \delta} M_j \in K^0_{ap}$. \newline
[Why?  By 1.15(2)].
\endproclaim
\bigskip

\proclaim{1.23 Claim}  Let $\bar N^\zeta,N_\zeta$ (for $\zeta \le \lambda)$
be as in 1.10.  Let $E \subseteq \lambda^+$ be a thin enough club of
$\lambda^+,\{ \varepsilon(\alpha):\alpha < \lambda^+\}$ enumerate
$\{ 0\} \cup E,H$ a 1-to-1 map from $N_\lambda$ onto $\lambda^+$ mapping
$N^\lambda_{\varepsilon(\alpha)}$ onto $\lambda \times \alpha$.  Let
$N^*_\alpha = H(N^\lambda_{\varepsilon(\alpha)}),N^* = \dsize \bigcup_{\alpha
< \lambda} N^*_\alpha$. \newline
1)  If $M \in K^0_{ap}$ then there is a lawful $f$ (see \cite{Sh:457},4.1)
which is an elementary embedding of $M$ into $N^*$ such that for
$\alpha \in \text{ Dom}(M),f(M \restriction (\lambda \times \alpha))
<^2_{f(M \restriction (\lambda \times \alpha + \lambda))} N^* \restriction
(\lambda \times \alpha)$.
\endproclaim
\bigskip

\demo{Proof}  Straightforward. \newline
[Saharon:  put old proof of 1.27 from AP here?
\enddemo
\bigskip

\noindent
But we want more, not only universality but also homogeneity.
\definition{1.24 Definition}  $K^1_{ap} = K^1_{ap}[T,\lambda]$ is the set of
$M \in K^0_{ap}$ such that for every $\alpha \in (0,\lambda^+)$, if
$\neg(M \subseteq \lambda \times \alpha)$ then $(M \restriction (\lambda
\times \alpha),M) \in K^{pr}_1$.

Let $\le_{K^1_{ap}}$ be $\le_{K^0_{ap}} \restriction K^1_{ap}$.
\enddefinition
\bigskip

\proclaim{1.25 Claim}  1) $\le_{K^1_{ap}}$ is a partial order on
$K^1_{ap}$. \newline
[Why?  By 1.22 and Definition 1.2]. \newline
2)  If $\langle M_i:i < \delta \rangle$ is $\le_{K^1_{ap}}$-increasing,
$\dsize \sum_{1 < \delta} \| M_i \| < \lambda$ \underbar{then} 
$M_i \le_{K^1_{ap}} \dsize
\bigcup_{j < \delta} M_j \in K^1_{ap}$. \newline
[Why?  By 1.23(1) and 1.19(1)].
\endproclaim
\bigskip

\proclaim{1.26 Claim}  Let $\bar N^\zeta,N_\zeta$ (for $\zeta \le \lambda)$
be as in 1.10.  Let $E \subseteq \lambda^+$ be a thin enough club of
$\lambda^+,\{ \varepsilon(\alpha):\alpha < \lambda^+\}$ enumerate
$\{ 0\} \cup E,H$ a 1-to-1 map from $N_\lambda$ onto $\lambda^+$ mapping
$N^\lambda_{\varepsilon(\alpha)}$ onto $\lambda \times \alpha$.  Let
$N^*_\alpha = H(N^\lambda_{\varepsilon(\alpha)}),N^* = \dsize \bigcup_{\alpha
< \lambda} N^*_\alpha$. \newline
1)  If $M \in K^1_{ap}$ then there is a lawful $f$ (see \cite{Sh:457},4.1)
which is an elementary embedding of $M$ into $N^*$ such that for
$\alpha \in \text{ Dom}(M),f(M \restriction (\lambda \times \alpha))
<^2_{f(M \restriction (\lambda \times \alpha + \lambda))} N^* \restriction
(\lambda \times \alpha)$.
\newline
2)  If $M_0 \le_{K^1_{ap}} M_1$ and $(M_0,f_0)$ is as in part (1)
\underbar{then}  we can find $f_1,f_0 \subseteq f_1$ such that $(M_1,f_1)$
is as in part (1).  Moreover, if $f_0 \cup (f_1 \restriction (M_1 \restriction
(M_1 \restriction (\lambda \times \alpha))$ has been determined we can
continue.
\endproclaim
\bigskip

\proclaim{1.27 Amalgamation Claim}  Assume $M_0 \le_{K^0_{ap}} M_\ell$ for
$\ell = 1,2$ and (for simplicity) $|M_1| \cap |M_2| = |M_0|$.
\underbar{Then} there is $M \in K^0_{ap}$ such that $M_1 \le_{K^0_{ap}} M$
and $M_2 \le_{K_{ap}} M$.
\endproclaim
\bigskip

\demo{Proof}  Follows from 1.26(1) + (2) ($Q$: domain?)
\enddemo
\bigskip

\proclaim{1.28 Claim}  $(K_{ap},\le^*)$ is a $\lambda$-system (see
\cite{Sh:457},\S4).  
\endproclaim
\bigskip

\demo{Proof}  Check.
\enddemo
\bigskip

\proclaim{1.29 Claim}  $(K_{ap},\le^*)$ is simple (see \cite{Sh:457},\S4).
\endproclaim
\bigskip

\demo{Proof}  Included in the proof of amalgamation (see last clause of
1.26(2)).
\enddemo
\bigskip

\proclaim{1.30 Claim}  If $M$ is a model of $T$ of cardinality $\lambda^+$
\underbar{then} for some $\Gamma \in K^{md}_\lambda,M$ can be elementarily
embedded into $M_\Gamma$.
\endproclaim
\bigskip

\demo{Proof}  Use 1.10 with $M = \dsize \bigcup_{\alpha < \lambda} M_\alpha$
so we get $N^*,N^*_\alpha(\alpha < \lambda)$ as in 1.18.  Check.
\enddemo
\bigskip

\demo{1.31 Proof of 1.2}  Use the above claims.
\enddemo
\newpage

\head {\S2 on the strong order properties and finitary versions} \endhead
\bigskip

\demo{2.0 Discussion}  By \cite{Sh:457}, for some non-simple (first order
complete) the answer to the following is yes:
\medskip
\roster
\item  "{$\bigoplus_T$}"  if $\lambda = \lambda^{< \lambda} > |T|,2^\lambda =
\lambda^+$, is there a ($\lambda$-complete), \newline
$\lambda^+$-c.c. forcing notion $Q$, $\Vdash_Q$ 
``univ$(\lambda^+,T) \le \lambda^{++} < 2^\lambda$"?
\item "{$\bigoplus'_T$}"  and by [?] if $\lambda = \lambda^{< \lambda} > |T|,
2^\lambda = \lambda^+$ is there a $\lambda$-complete $\lambda^+$-c.c. forcing
notion $Q,\Vdash_Q$ ``univ$(\lambda^+,T) = 1, \lambda^+ < 2^\lambda$"?
\endroster
\medskip

We know that for theories $T$ with the strict order property the answer is no
(by \cite{KjSh:409}, or see \cite{Sh:457},\S3).  We would like to characterize
the answer by a natural property of $T$ (hence show that the answer to all
reasonable variants is the same, e.g. does not depend on $\lambda,
\oplus_T \equiv \oplus'_T$, etc.)  So the results we mention above give a
lower bound (simple theories $+T_{qef} + T_{trf}$) and an upper bound
(failure of the strict order property) to the family of $T$'s with a positive
answer.  However, we can lower the upper bound.  We suggest below a strictly
weaker property.  From another point of view, a major theme of cite{Sh:a},
\cite{Sh:c} was to find natural dividing lines for the family of first
order theories (so the main ones there were stable, superstable and also
NTOP, deepness NOTOP).  Now \cite{Sh:93} suggests another one: simplicity.
Note that the negation of simple, the tree property has been touched upon
in \cite{Sh:a} but there were conclusions only for one side.  \cite{Sh:93}
establishes this dividing line by having consequences for both the property
and its negation and having ``semantical characterization" for $T$ simple:
when $|T| \le \kappa < \lambda = \lambda^{< \lambda} < \mu = \mu^\kappa$ we
can force by a $\lambda^+$-c.c. $\lambda$-complete forcing notion $Q$ that
$2^\lambda > \mu$ and every model of $T$ of cardinality $\mu$ can be
extended to a $\kappa^+$-saturated one, and the tree property implies a
strong negation.  Of course, both the inner theory and such ``outside",
``semantical" characterization are much weaker than those for stable
theories.

The strict order property has no such results only several consequences.
We suggest below weaker properties (first the strong order property then
the $n$-version of it for $n < \omega$) which has similar consequences and
so may be the right dividing line (for some questions).  Remember (this is
in equivalent formulations).
\enddemo
\bigskip

\definition{2.1 Definition}  $T$ has the strict order property if some
formula $\varphi(\bar x,\bar y)$ (with $\ell g \bar x = \ell g \bar y)$
define in some model $M$ of $T$, a partial order with infinite chains.
\enddefinition
\bigskip

\definition{2.2 Definition}  1) A first order complete $T$ has the strong
order property if some sequence $\bar \varphi = \langle \varphi_n(\bar x^n;
\bar y^n):n < \omega \rangle$ of formulas exemplifies it which means that
for every $\lambda$:
\medskip
\roster
\item "{$(*)^\lambda_{\bar \varphi}(a)$}"  $\ell g \bar x^n = \ell g \bar y^n$
are finite, $\bar x^n$ an initial segment of $\bar x^{n + 1}$
$$
\bar y^n \text{ an initial segment of } \bar y^{n + 1}
$$
\item "{$(b)$}"  $T,\varphi_{n+1}(\bar x^{n+1},\bar y^{n+1}) \vdash
\varphi_n(\bar x^n,\bar y^n)$
\item "{$(c)$}"  for $m \le n, \neg(\exists \bar x^{n,0}\cdots \bar x^{n,m-1})
[\bigwedge \{ \varphi_n(\bar x^{n,k},\bar x^{n,\ell}):k = \ell + 1
\text{ mod } m\}]$ belongs to $T$
\item "{$(d)$}"  there is a model $M$ of $T$ and $\bar a^n_\alpha \in M$ (of
length $\bar y^n$, for $n < \omega,\alpha < \lambda$) such that
$\bar a^n_\alpha = \bar a^{n+1}_\alpha \restriction \ell g \bar y^n$ and
$M \models \varphi_n[\bar a^n_\beta,\bar a^n_\alpha]$ for $n < \omega$ and
$\alpha < \beta < \lambda$.
\endroster
\medskip

\noindent
2)  The finitary strong order property is defined similarly but $\bar x^n =
\bar x,\bar y^n = \bar y^n$. \newline
3)  We use the shorthand SOP, FSOP and for the ngation NSOP, NFSOP (similarly
later for NSOPn).
\enddefinition
\bigskip

\proclaim{2.3 Claim}  1) The strict order property implies finitary strong
order property which implies the strong order property. \newline
2)  There is a first order complete $T$, which has the strong order property
(even the finitary one) but not the strict order property. \newline
3) Also some first order complete $T$ has the strong order property but not
the finitary strong order property, i.e. no $\langle \varphi_n(\bar x,\bar y):
n < \omega \rangle$ exemplifies it (i.e. with $\ell g \bar x_n$ constant).
\endproclaim
\bigskip

\demo{Proof}  1) Immediate. \newline
2)  For $\ell \le n < \omega$ let $<_{n,\ell}$ be a two-place relation.
Let $<_n = <_{n,0}$.  Let  $T_0$ say:
\medskip
\roster
\item "{(a)}"  $x <_{n,m-1} y \Rightarrow x <_{n,m} y$
\item "{(b)}"  $x <_{n,n} y$
\item "{(c)}"  $\neg(x <_{n,n-1} x)$
\item "{(d)}"  if $\ell + k + 1 = m \le n$ then $x <_{n,\ell} y \and
y <_{n,k} z \Rightarrow x <_{n,m} z$.
\endroster
\medskip

\noindent
We shall now prove that $T_0$ has the amalgamation property; it also has the
point embedding property (as the latter is easier we leave its checking to
the reader).

Now suppose $M_i \models T_0,M_0 \subseteq M_i$ for $i = 0,1,2$ and
$M_1 \cap M_2 = M_0$.  We define a model $M$: its universe is $M_1 \cup
M_2$ and

$$
\align
<^M_{n,m} = \biggl\{(a,b) \in M \times M:&\text{if }m < n \text{ then for
some }i \in \{ 1,2\} \text{ we have}: \\
  &(a,b) \in <^{M_i}_{n,m} \text{ or } a \in M_i \backslash M_0,b \in
M_{3-i} \backslash M_0 \\
  &\text{and for some } c \in M_0 \text{ and } \ell,k 
\text{ we have}: \\
  &m = \ell + k + 1,(a,c) \in <^{M_i}_{n,\ell},(c,b) \in <^{M_{z-i}}_{n,k} 
\biggr\}.
\endalign
$$
\medskip

\noindent
Now clearly $M$ extends $M_1$ and $M_2$: trivially $<^M_{n,m} \restriction
M_i = <^{M_i}_{n,m}$.  Is $M$ a model of $T_0$?  Let us check. 
\enddemo
\bigskip

\noindent
\underbar{Clause $(a)$ holds}:  For $x,y \in M_i$ as $M_i \subseteq M$; for
$i = 1,2$ and $x \in M_i \backslash M_0,y \in M_{3-\ell} \backslash M_0$,
without loss of generality $m < n$; let $c \in M$ witness $(a,b) \in
\le^M_{n,m-1}$ i.e. for some $\ell,k$ we have $\ell + k + 1 = m-1,(a,c) \in
<^{M_i}_{n,\ell}$ and $(c,b) \in <^{3-i}_{n,k}$.  Now by clause (a) applied
to $M_i,(a,c) \in <^{M_i}_{n,\ell + 1}$ now apply the definition to get
$(a,b) \in <^M_{n,(\ell + 1)+k+1} = <^M_{n,m}$.
\bigskip

\noindent
\underbar{Clause $(b)$ holds}:  Check as defining $<^M_{n,m}$ we say: ``if
$m < n$ then ..." so if $n = m$ there is no requirement.
\bigskip

\noindent
\underbar{Clause $(c)$}:  As $M_i \subseteq M$ and $M_i \models T_0$.
\bigskip

\noindent
\underbar{Clause $(d)$}:  Check by cases, i.e. for some $i \in \{1,2\}$ one
of the following cases hold.
\medskip
\roster
\item  $\{x,y,z\} \subseteq M_i$: \newline
use ``$M_i$ is a model of $T_0$ and $M_i$ a submodel of $M$".
\item  $\{x,y\} \subseteq M_i,\{y,z\} \subseteq M_{3-i}$: \newline
use the definition of $<^M_{n,m}$.
\item  $y \in M_i \backslash M_0,\{x,z\} \subseteq M_{3-i} \backslash
M_0$.
\endroster
\medskip

\noindent
As $x <_{n,\ell} y$ there are $\ell_1,\ell_2$ and $x_1 \in M_0$ such that
$x <_{n,\ell_1} x_1,x_1 \le_{n,\ell_2} y$ and $\ell_1 + \ell_2 + 1 = \ell$.

As $y <_{n,k} z$ there are $k_1,k_2$ and $z_1 \in M_0$ such that
$y \le^{M_i}_{n,k_1} z_1, z_1 \le^{M_{3-i}}_{n,k_2} z,k_1 + k_2 + 1 = k$.
In $M_i$ we have $x_1 <^{M_i}_{n,\ell_2} y <^{M_i}_{n,k_1} z_1$ hence
$x_1 \le^{M_i}_{n,\ell_2 + k_1 + 1} z_1$ and as $\{ x_1,z_1\} \subseteq M_0
\subseteq M_i$ clearly $x_1 <^{M_0}_{n,\ell_2 + k_1 +1} z_1$.  Now in
$M_{3-i}$ we have $x <^{M_{3-i}}_{n,\ell_1} x_1 <^{M_{3-i}}_{n,\ell_1} x_1
<^{M_{3-i}}_{n,\ell_2 + k_1 + 1} z_1$ hence
$x \le^{M_{3-i}}_{n,\ell_1 + \ell_2 + k_1 + 2} z_1$ so 
$x <^{M_{3-\ell}}_{n,\ell_1 + \ell_2 + k_1 + 2} z_1 <^{M_{3-\ell}}_{n,k_2}z$
hence $x <^{M_{3-\ell}}_{n,\ell_1 + \ell_2 + k_1 + k_2 + 3}z$ but
$\ell_1 + \ell_2 + k_1 + k_2 + 3 = \ell + k + 1 = m$ so
$x <^{M_{3-\ell}}_{n,m}z$ as required.
\medskip
\roster
\item "{(4)}"  $y \in M_i \backslash M_0, x \in M_{3-i} \backslash M_0,
z \in M_0$. \newline
Similar to case (3) but with no $x_1$.
\item "{(5)}"  $y \in M_i \backslash M_0,x \in M_0,z \in M_{3-i} \backslash
M_i$. \newline
Similar to case (3) but with no $z_1$.
\endroster
\medskip

Let $T$ be the model completion of $T^0$; easy to check that it exists and
has elimination of quantifilters.  Let $\varphi_n(x,y) = \dsize
\bigwedge_{\ell \le n} x <_\ell y$ (remember $x <_\ell y$ means $x
<_{\ell,\delta} y$) now $\langle \varphi_n:n < \omega \rangle$ exemplifies
that $T$ has the (finitary) strong order property.  On the other hand we
shall show that for every $n(*) < \omega$ the theory $T_{n(*)} =: T
\restriction \{ <_{n,\ell}:\ell \le n \le n(*)\}$ does not have the strict
order property (as $T = \dsize \bigcup_{n < \omega} T_n$, this clearly 
implies that $T$ does not have the strict order property).  First note that
also $_{n(*)}$ has elimination of quantifiers and then check directly.
\bigskip

(3) Let $T^0$ say:
\medskip
\roster
\item "{(a)}"  $P_n$ (for $n < \omega$) are pairwise disjoint ($P_n$
unary predicates)
\item "{(b)}"  $F_n$ a partial one place function from $P_{n+1}$ into
$P_n$
\item "{(c)}"  $<_{n,\ell}$ are two-place relations on $P_n$ for $\ell \le
n < \omega$; and let $<_n = <_{n,0}$
{\roster
\itemitem{ $(\alpha)$ }  $x <_{n,m-1} y \rightarrow x <_{n,m} y$
\itemitem{ $(\beta)$ }  $P_n(x) \and P_n(y) \rightarrow x <_{n,n} y$
\itemitem{ $(\gamma)$ }  $\neg(x <_{n,n-1} x)$
\itemitem{ $(\delta)$ }  if $\ell + k + 1 = m \le n$ then: $x <_{n,\ell} y
\and y <_{n,k} z \rightarrow x <_{n,m} z$
\endroster}
\item "{(d)}"  $x <_{n+1,\ell} y \rightarrow F_n(x) <_{n,\ell} F_n(y)$.
\endroster
\medskip

Again $T$ will be the model completion of $T^0$ and it has elimination of
quantifiers and we shall use $\bar x_n = \langle x_i:i < n \rangle,
\bar y_n = \langle y_i:i < n \rangle$ and $\varphi_n(\bar x_n,\bar y_n) =
\dsize \bigwedge_{i < n} F_i(x_{i+1}) = x_i \and \dsize \bigwedge
_{i < n} F_i(y_{i+1}) = y_i \and \dsize \bigwedge_{i < n} x_i <_i y_i$.
\hfill$\square_{2.3}$
\bigskip

\proclaim{2.4 Claim}  1) The following are equivalent (for $\lambda \ge |T|$):
\medskip
\roster
\item "{$(A)$}"  $T$ has the strong order property
\item "{$(B)_\lambda$}"  there is a $\lambda^+$-saturated model $M$ of $T$,
a $L_{\infty,\lambda^+}$-formula $\varphi = \varphi(\bar x,\bar y),
\varepsilon = \ell g \bar x = \ell g \bar y \le \lambda$, possible with
$\le \lambda$ parameters, such that in $M,\varphi$ defines a partial linear
order with a chain of length $\ge \beth_2(\lambda)^+$. 
\endroster
\medskip

\noindent
2)  The following are equivalent $(\lambda \ge |T|$):
\medskip
\roster
\item "{$(A)'$}"  $T$ has the finitary strong order property
\item "{$(B)'_\lambda$}"  like $(B)_\lambda$ but $\varepsilon < \omega$.
\endroster
\endproclaim    
\bigskip

\demo{Proof}  1) \underbar{$(A) \Rightarrow (B)_\lambda$}  Straight: for a
given $\bar \varphi = \langle \varphi_n(\bar x_n,\bar y_n):n < \omega
\rangle$, let $\bar x,\bar y$ be the limit of $\bar x_n,\bar y_n$
respectively and write $\psi^*(\bar x,\bar y) = \dsize \bigvee_m (\exists
\bar z_0,\dotsc,\bar z_m)[\bar x = \bar z_0 \and \bar y = \bar z_m \and
\dsize \bigwedge_{\ell < m} \varphi_\omega(\bar z_\ell,\bar z_{\ell + 1}]$
where $\varphi_\omega(\bar x,\bar y) =: \dsize \bigwedge_n \varphi_n(\bar x,
\bar y)$.
\bigskip

\noindent
\underbar{$(B)_\lambda \Rightarrow (A)$}  Let $\bar a_\alpha \in 
{}^\varepsilon M$ for $\alpha < \beth_2(\lambda)^+$ form a chain.  Without
loss of generality the order $\varphi$ defines is strict (i.e. $\vdash \varphi
(\bar x,\bar x)$) and no parameters (just add them to the $\bar a_\alpha$'s).
By Erdos Rado theorem without loss of generality for some type $q = q(\bar x,
\bar y)$ for all $\alpha < \beta < \omega$ the sequence $\bar a_\alpha 
\char 94 \bar a_\beta$ realizes $q$.

For every $n,\bigcup \{ q(\bar x_\ell,\bar x_k):k = \ell + 1 \text{ mod } n
\text{ and } k,\ell < n\}$ cannot be realized in $M$ (as if $\bar b_0 \char 94
\cdots \char 94 \bar b_{n-1}$ realizes if we get a contradiction to
``$\varphi(\bar x,\bar y)$ defines a strict partial order").  By saturation
there is $\varphi^0_n(\bar x,\bar y) \in q(\bar x,\bar y)$ such that
$\{ \varphi^0_n(\bar x_\ell,\bar x_k):k = \ell + 1 \text{ mod } n$ and
$k,\ell < n\}$ is not realized in $M$.  The rest should be clear.
\newline
2)  Left to the reader. \hfill$\square_{2.4}$
\enddemo
\bigskip

\definition{2.5 Definition}  1) $T$ has the $n$-stronger order property
(SOP$_n$) \underbar{if} there is a formula $\varphi(\bar x,\bar y)$ having
this property for $T$ which means: $\ell g \bar x = \ell g \bar y$ (allowing
parameters changes nothing) and there is a model $M$ of $T$ and
$\bar a_k \in {}^{\ell g \bar x}M$ for $k < \omega$ such that:
\medskip
\roster
\item "{(a)}"  $M \models \varphi[\bar a_k,\bar a_m]$ for $k < m < \omega$
\item "{(b)}"  $M \models \neg \exists \bar x_0 \cdots \bar x_{n-1}
(\bigwedge \{ \varphi(\bar x_\ell,\bar x_k):\ell,k < n \text{ and }
k = \ell + 1 \text{ mod } n\})$.
\endroster
\medskip

\noindent
2) ``$T,\varphi(\bar x,\bar y)$ have the SOP$_{\le n}$" is defined similarly
except that in $(b)$ we replace $n$ by each $m \le n$.
\enddefinition
\bigskip

\proclaim{2.6 Claim}  SOP $\Rightarrow$ SOP$_{n+1}$, SOP$_{n+1} \Rightarrow$
SOP$_n$, SOP$_{\le n+1} \Rightarrow$ SOP$_{\le n}$ \newline
and
SOP$_n \Leftrightarrow$ SOP$_{\le n}$ for any given $T$, (we did not say
``for any $\varphi$").
\endproclaim
\bigskip

\demo{Proof}  The first clause is immediate.  The second clause is straight
too: \newline
let $\varphi(\bar x,\bar y),M,\langle \bar a_m:m < \omega \rangle$
exemplify SOP$_{n+1}$ and without loss of generality the sequence
$\langle \bar a_m:m < \omega \rangle$ is an indiscernible sequence.  Does
$M \models (\exists \bar x_0,\dotsc,\bar x_{n-1})[\bar x_0 = \bar a_1 \and
\bar x_{n-1} = \bar a_0 \and \dsize \bigwedge \{ \varphi(\bar x_\ell,
\bar x_k):\ell,k < n \text{ and } k = \ell + 1 \text{ mod } n\})$?  If the
answer is yes we can replace $\bar a_1$ by $\bar a_2$ (by indiscernability),
let $\bar c_0,\dotsc,\bar c_{n-1}$ be as required above on
$\bar x_0,\dotsc,\bar x_{n-1}$ and $\bar b_0 =: \bar a_1,\bar b_1 =: \bar a_2
(= \bar c_0),\bar b_2 =: \bar c_1,\dotsc,\bar b_{n-1} =: \bar c_{n-2}$,
$\bar b_n =: \bar c_{n-1} = a_0$; now they satisfy the requirement mentioned
in $(b)$ of 2.5(1) on $\bar x_0,\dotsc,x_n$ (for SOP$_{n+1}$), contradicting
clause $(b)$ of 2.5(1).  So assume ``no" and now $\varphi'(\bar x,\bar y)$
have SOP$_n$ for $T$ where: $\varphi'(\bar x,\bar y) =: \varphi(\bar x,\bar y)
\and \neg (\exists \bar x_0,\dotsc,x_{n-1})[\bar x_0 = \bar x \and \bar x_1
= \bar y \and \dsize \bigwedge \{ \varphi(\bar x_\ell,\bar x_k):\ell,k < n$
and $k = \ell \text{ mod } n\}]$.
\medskip

As for SOP$_n \Leftrightarrow SOP_{\le n}$, the implications $\rightarrow$ is
really included in the proof above, (i.e. by it, if $\langle \bar a_\ell:
\ell < \omega \rangle,\varphi_n$ exemplifies SOP$_n$, for some $\varphi_{n-1}$
we have $\langle \bar a_\ell:\ell < \omega \rangle,\varphi_{n-1}$ exemplifies
SOP$_{n-1}$ (with $n,n-1$ here corresponding to $n+1,n$ there), and we can
define $\varphi_{n-2},\cdots$ similarly; now $\langle \bar a_\ell:\ell < 
\omega
\rangle, \dsize \bigwedge_{i \le n} \varphi_i$ exemplifies SOP$_{\le n}$.
The implication $\Leftarrow$ is trivial.  Now the third clause
SOP$_{\le n+1} \Rightarrow$ SOP$_{\le n}$ is trivial (read the definition).
\hfill$\square_{2.6}$
\enddemo
\bigskip

\proclaim{2.7 Claim}  Let $T$ be complete.  If $T$ has SOP$_3$ then $T$ has
the tree property (i.e. is not simple).
\endproclaim
\bigskip

\demo{Proof}  Let $\kappa = \text{ cf}(\kappa) > |T|$ and $\lambda > \kappa$
be a strong limit singular cardinal of cofinality $\kappa$.  Let
$J = {}^\kappa \lambda,I = \{ \eta \in {}^\kappa \lambda:\eta(i) = 0
\text{ for every } i < \kappa \text{ large enough}\}$.  Let $\varphi
(\bar x,\bar y)$ exemplify the SOP$_3$.  By the definition we can find a
model $M$ of $T$ and $\bar a_\eta \in M$ (for $\eta \in J$) such that:
\medskip
\roster
\item "{$(*)$}"  $\eta <_{\ell x} \nu$ in $I \Rightarrow M \models \varphi
[\bar a_\eta,\bar a_\nu]$.
\endroster
\medskip

Without loss of generality $\| M \| \ge \lambda,M$ is $\kappa^+$-saturated.
So for every \newline
$\eta \in {}^\kappa(\lambda \backslash \{ 0 \}) \backslash I$ we
can find $\bar a_\eta \in M$ such that it realizes \newline
$p_\eta = \{ \varphi
(a_{(\eta \restriction i) \char 94 0_{[i,\kappa)}},\bar x) \and \varphi
(\bar x,\bar a_{(\eta \restriction i) \char 94 \langle \eta(i)+1 \rangle
\char 94 0_{[i,\kappa)}}):i < \kappa\}$.  But if \newline
$\eta_1 <_{\ell x} \eta_2
\in {}^\kappa(\lambda \backslash \{0\}$) then we can find $\nu,\rho \in I$
such that: \newline
$\eta_1 <_{\ell x} \nu <_{\ell x} \rho <_{\ell x} \eta_2$ and
$\varphi(\bar x,\bar a_\nu) \in p_{\eta_1},\varphi(\bar a_\rho,\bar x) \in
p_{\eta_2}$ and by $(*)$ we have $M \models \varphi[\bar a_\nu,\bar a_\rho]$,
so $p_{\eta_1} \cup p_{\eta_2}$ is contradictory (by clause $(b)$ of 2.5(1)
for ``$\varphi$ have the SOP$_3$").  So $\langle p_\eta:\eta \in {}^\kappa
(\lambda \backslash \{ 0 \})$ are pairwise contradictory, $|p_\eta| = \kappa$,
and $\lambda^\kappa > \lambda = \lambda^{< \kappa} > 2^{|T|}$ and
$\dsize \bigcup\{\text{Dom }p_\eta:\eta \in {}^\kappa(\lambda \backslash
\{ 0\})\}$ has cardinality $\le \lambda$ and $\kappa > |T|$.

By \cite{Sh:a},III,7.7 = \cite{Sh:c},III,7.7,p.141 this implies that $T$ has
the tree \newline
property. \hfill$\square_{2.7}$
\enddemo
\bigskip

\proclaim{2.8 Claim}  1) The theory $T_n =: T \restriction \{ <_{n,\ell}:\ell
\le n\}$ from 2.3(2) has SOP$_n$ but not SOP$_{n+1}$. \newline
2)  $T^{\text{mc}}_{\text{trf}}$, the model completion of the theory of
triangle free graphs has SOP$_3$ but not SOP$_4$. \newline
3) For $n \ge 3$ the model completion $T^{\text{mc}}_n = T^{\text{mc}}
_{\text{dcf}(n)}$ of the theory $T_n = T_{\text{dcf}(n)}$ of graphs 
(= directed graphs, no loops or multiple edge for simplicity) with no
directed circle of length $\le n$ has SOP$_n$ but not SOP$_{n+1}$. \newline
4) For odd $n \ge 3$, the model completion $T^{\text{mc}}_n =
T^{\text{mc}}_{\text{ocf}(n)}$ of the theory $T_n = T_{\text{ocf}(n)}$ of
graphs with no \underbar{odd} circle of length $\le n$, has SOP$_n$ but not
SOP$_{n+1}$. \newline
5) For $n \ge 3$, the model completion $T^{\text{mc}}_{\text{cf}(n)}$ of the
theory $T_n = T_{\text{cf}(n)}$ of graphs with no circles of length $\le n$,
has SOP$_3$ but not SOP$_4$. \newline
6) The theory $T_{\text{qcf}}$ (see \cite{Sh:457} does not have SOP$_3$ (but
is not simple).
\endproclaim
\bigskip

\remark{2.8A Remark}  1) Note that univ$(\lambda,T^{\text{mc}}_{\text{cf(n)}})
= \text{univ}(\lambda,T_{\text{cf(n)}})$. \newline
2) For those theories, $D(T^{\text{mc}})$ is an uncountable; they have no
universal model in $\lambda < 2^{\aleph_0}$.
\endremark
\bigskip

\demo{Proof}  1) Proved really in 2.3. \newline
2) This is included in part (5). \newline
3), 4), 5)  We discuss the existence of model completion later; note that the
meaning of $T_n$ depends on the part we are proving. \newline
Let $xRy$ mean $(x,y)$ is an edge; when we say $(x,y)$ is an edge, for
graphs we mean $\{x,y\}$ is an edge.  Let $\bar y = \langle y_\ell:\ell < n
\rangle,\varphi(\bar x,\bar y) = \dsize \bigwedge_{\ell < n-1} x_\ell
Ry_{\ell + 1} \and x_{n-1} Ry_0$.  First we note there $T_n \vdash \neg
(\exists \bar x_0,\dotsc,\bar x_{n-1}) \dsize \bigwedge \{\varphi(\bar x_i,
\bar x_k):\ell,k < n,k = \ell + 1 \text{ mod } n\}$, otherwise there are
$M \models T_n$ and $\bar a_\ell = \langle a_{\ell,0},\dotsc,a_{\ell,n-1}
\rangle \in {}^n M$ as forbidden but then $a_{0,0},a_{1,1},\dotsc,a_{n-1,n-1}$
is a circle, so in all cases this is impossible.
\medskip

For parts 3), 4) let $M$ be the following model of $T_n$; elements \newline
$a^\ell_i \, (i < \omega,\ell < n), R = \{ (a^\ell_i,a^{\ell + 1}_j):
i < j < \omega,\ell < n - 1\} \cup \{(a^{n-1}_i,a^0_j):i < j < \omega\}$ (but
for graphs we put all such pairs and the inverted pair as $R$ should be
symmetric and irreflexive relation).  Clearly for $i < j < \omega, M \models
\varphi[\bar a_i,\bar a_j]$ where $\bar a_i = \langle a^0_i,\dotsc,a^{n+1}
_i \rangle$.  Lastly $M \models T$: for part (3) as $R$ is not symmetric the
absence of any circle should be clear, $M \models a^{\ell(1)}_{i(1)}
Ra^{\ell(2)}_{i(1)} \Rightarrow i(1) < i(2)$; for part (4) there are circles
but even or long and 
$M \models a^{\ell(1)}_{i(1)}Ra^{\ell(2)}_{i(2)} \Rightarrow
\ell(1) = \ell(2) + 1 \text{ mod }n$, so $T_{\text{dcf}(n)},T_{\text{ocf}(n)}$
(and $T^{\text{mc}}_{\text{ocf}(n)}$) has
even or long $M \models a^{\ell(1)}_{i(1)}Ra^{\ell(2)}_{i(2)} \Rightarrow
\ell(1) = \ell(2) + 1 \text{ mod }n$.  So $T_{\text{dcf}(n)},
T_{\text{ocf}(n)}$ (and $T^{\text{mc}}_{\text{dcf}(n)}$ and
$T^{\text{mc}}_{\text{ocf}(n)}$) has SOF$_n$.
\medskip

Let $n = 3$.  Now $T_{\text{cf}(3)} = T_{\text{ocf}(3)}$ so we can ignore
part (5).  Also $T_{\text{dcf}(n)},T_{\text{ocf}(n)}$ has the amalgamation
property and joing embedding property.  Thus, it is enough to show that
$T^{\text{mc}}_{\text{dcf}(n)},T^{\text{mc}}_{\text{ocf}(n)}$ fails the
SOP$_4$.  As $T^{\text{md}}$ has elimination of quantifiers the reader can
check directly that $T_n$ does not have SOP$_4$.
\bigskip

Let $n > 3$.  
Though $T^{\text{mc}}_n$ does not have elimination of quantifiers, every 
formula is equivalent to a Boolean combination of formulas of the form:
$x = y,xRy$ for $m < n,\varphi_m(x,y) =: (\exists x_0,\dotsc,x_m)[x = x_0
\and y = x_m \bigwedge \dsize \bigwedge_{\ell < m} x_\ell Rx_{\ell + 1}]$
(i.e. the distance from $x$ to $y$ is $\le m$, directed from $x$ to $y$ in
the case of di-graphs). For part (5) of 2.8, we should add for $\ell < m 
< n/2,\ell > 0$ a partial function $F_{m,\ell}$ defined by:
$F_{m,\ell}(x,y) = z$ iff there are $t_0,\dotsc,t_m$ with no repetition such
that $x = t_0,y = t_m,z = t_\ell$ and $\dsize \bigwedge_{\ell < m}
t_\ell Rt_{\ell + 1}$ and lastly $\psi_{m,\ell}(x,y) =: (\exists z)
[F_{m,\ell}(x,z) = z]$.  Let $T^2_n$ be the set of obvious (universal) axioms
for those relations.  Then easily $T^2_n$ has amalgamation and has model
completion, $T^1_n$ which has elimination of quantifiers (but the closure
of a finite set under those functions may be infinite).  Moreover, assume
$M \models T^2_n,\langle \bar a_m:m < \omega \rangle$ is an indiscernible
sequence in $M,\bar a_m = \langle a^m_\ell:\ell < k \rangle$, with
$k < \omega$.  Then there is $w \subseteq k$ such that $[a^m_\ell =
a^{m+1}_\ell \Leftrightarrow \ell \in \omega]$ and without loss of generality
$[\ell_1 < \ell_2 \Rightarrow a^m_{\ell_1} \ne a^m_{\ell_2}]$.  Let for
$u \subseteq \omega,M_u$ be the submodel of $\mu$ generated by
$\dsize \bigcup_{m \in u} \bar a_m$ for parts (3), (4), $M_u = \bigcup
\{ M_v:v \subseteq u \text{ and } |v| \le 1\}$ so things are simple.
By the indiscernibility (increasing the $\bar a_m$'s e.g. taking $\omega$
blocks) without loss of generality
\medskip
\roster
\item "{$(*)$}"  $M_u \cap M_\nu = M_{u \cap \nu}$ and the universe of
$M_{\{m\}}$ is the range of $\bar a_m$.
\endroster
\medskip

Let $m = n$ for parts (3), (4) of 2.8, $m = 3$ for part (5).  For part (5)
note: the distance between $a^0_{\ell_0},a^1_{\ell_1}$ is $> {\frac n4}$.
\newline
[Why?  If not there is a path $C^{i,j}$ of length $\le {\frac n4}$ for
$a^i_{\ell_0}$ to $a^j_{\ell_1}$, now $C^{0,3} \cup C^{1,2} \cup C^{1,4} \cup
C^{0,4}$ is a circle of length $\le n$, may cross itself but still there is
a too small circle]. \newline

We can now define models $N_{\{ \ell\}}$ (for $\ell < n + 1),
N_{\{\ell,\ell + 1\}}(\ell < n)$ and $N_{\{n,0\}}$ and isomorphisms
$h_\ell,g_\ell(\ell < n+1)$ such that:
\medskip
\roster
\item "{(a)}"  for $\ell < n+1 h_\ell$ an isomorphism from $M_{\{ \ell\}}$
onto $N_{\{\ell\}}$
\item "{(b)}"  for $\ell < n,g_\ell$ an isomorphism $M_{\{\ell,\ell + 1\}}$
onto $N_{\{\ell,\ell + 1\}}$ extending $h_\ell,h_{\ell + 1}$.
\item "{(c)}"  for $\ell = n,g_\ell$ an isomorphism from $M_{\{n,n + 1\}}$
onto $N_{\{n,0\}}$ extending $h_n$ and $h_0 \circ f$ where $f$ is the
isomorphism from $M_{\{n + 1\}}$ onto $M_{\{0\}}$ taking $\bar a_{n+1}$ onto
$\bar a_0$.
\item "{(d)}"  $N_\emptyset =: g_\ell(M_\emptyset)$ does not depend on
$\ell$.
\item "{(e)}"  $N_u \cap N_\nu = N_{u \cap \nu}$ if $u,\nu$ are among
$\emptyset,\{ \ell \},\{m,m+1\},\{n,0\}(\ell < n + 1,m < n)$.
\endroster
\medskip

\noindent
Now,
\medskip
\roster
\item "{$\bigotimes$}"  There is a model of $T^2_n$ extending all
$N_{\{\ell,\ell + 1\}},N_{\{n,0\}}(\ell < n)$.
\endroster
\medskip

This is enough for showing that $T^1_n$ lacks the SOP$_{m+1}$. \newline
Lastly the reader can check that $T^{\text{mc}}_{\text{cf}(n)}$ has
SOP$_3$ [choose $k \in ({\frac n3},n),\bar a_\ell = \langle a_{\ell,0}
\rangle,\ell_1 < \ell_2 \Rightarrow \varphi_k(a_{\ell,0},a_{\ell_2,0})$].
\hfill$\square_{2.8}$
\enddemo
\bigskip

\proclaim{2.9 Theorem}  Let $T$ be first order complete, $\lambda \ge |T|$
and $T$ has the SOP$_3$.  \underbar{Then}: \newline
1)  $T$ is maximal in the Keisler order $\triangleleft_\lambda$, i.e. for
a regular filter $D$ on $\lambda$ and some (= every) model $M \models T$ we
have $M^\lambda/D$ is $\lambda^+$-saturated \underbar{iff} $D$ is a good
ultrafilter. \newline
2) Moreover, in 2.10 $T$ is $\triangleleft^\ell$-maximal, (see Definition 2.10
below). \newline
We delay the proof.
\bigskip
\endproclaim
\bigskip

\remark{Remark}  The order $\triangleleft$ was introduced and investigated
by Keisler \cite{Ke76}; further investigated in \cite{Sh:42},
\cite{Sh:a},CH.VI, new version \cite{Sh:c},CH.VI.  The following is a
generalization. 
\endremark
\bigskip

\definition{2.10 Definition}  1) For models $M_0,M_1$ we say $M_0,
\triangleleft^*_\lambda M_1$ if the following holds: for some model
${\frak B}_0$ in which $M_0,M_1$ are intepreted (so $M_i = 
M^{{\frak B}_0}_i)$, for every elementary extension ${\frak B}$ of
${\frak B}_0$, which is $(\aleph_0 + | \tau(M_0)| + |\tau(M_1)|)^+$-saturated
we have: $[M^{\frak B}_1$ is $\lambda^+$-saturated $\Rightarrow
M^{\frak B}_0$ is $\lambda^+$-saturated]. \newline
2)  $M_0 \triangleleft^* M_1$ if for every $\lambda \ge \aleph_0 + |\tau
(M_0)| + |\tau(M_1)|$ we have $M_0 \triangleleft^*_\lambda M_1$. \newline
3)  Using the superscript $\ell$ instead of $*$ means in the saturation we
use only $\varphi$-types for some $\varphi = \varphi(\bar x,\bar y)$ (so any
$\varphi$ is O.K., but for each type $\varphi$ is constant) and omit the
saturation demand on ${\frak B}$. \newline
4)  For complete theories $T_1,T_2$ we say $T_1 \triangleleft^*_\lambda T_2$
\underbar{if} for every model $M_1$ of $T_1$ for some model $M_2$ of $T_2$,
$M_1 \triangleleft^*_\lambda M_2$.  Similarly for $T_1 \triangleleft^*_\lambda
T_2,T_1 \triangleleft^\ell_{(\lambda)} T_2$.
\enddefinition
\bigskip

\demo{2.11 Observation}  1) In 2.11(1) we can just use ${\frak B}_0$ of
the form \newline
$(H(\chi),\in,<^*_\chi,M_0,M_1)$ with $\chi$ strong limit. \newline
2) $\triangleleft^*_\lambda$ is a partial order, also $\triangleleft^\ell
_*,\triangleleft^\ell$ are partial orders; $M \triangleleft^*_\lambda M$, and
if $M_0$ is interpretable in $M_1$ then $M_0 \triangleleft^*_\lambda M_1$.
\newline
2A) For models of countable vocabulary, similar statements hold for
$\triangleleft^*$ (without the countability if $|\tau(M_1)| > |\tau(M_0)|
+ |\tau(M_2)| + \aleph_0$, we can get a silly situation). \footnote{ so to
overcome this, we may in Definition 2.10(2) replace ``every $\lambda >
\cdots$" by `` every large enough $\lambda$"} \newline
3) If $\lambda \ge \aleph_0 + |\tau(M_0)| + |\tau(M_1)|$ then:
$M_0 \triangleleft^*_\lambda M_1$ iff for every finite $\tau \subseteq
\tau(M_0),M_0 \restriction \tau \triangleleft^*_\lambda M_1$. \newline
4) $M_1 \otimes^\ell_\lambda M_2 \Rightarrow M_1 \triangleleft^*_\lambda
M_2$. \newline
5) Parallel results hold for theories. \newline
6) Any (complete first order) theory of any infinite linear order is
$\triangleleft^\ell$-maximal hence $\triangleleft^*_\lambda$-maximal for
every $\lambda \ge |T| + \aleph_0$. \newline
7) All countable stable theories without the f.c.p. (e.g. $T = Th(\omega_1
=)$) are $\triangleleft^\ell$-equivalent. \newline
8) All countable stable theories with the f.c.p. are equivalent \newline
(e.g.
$T_{\text{eq}} = Th(\dsize \bigcup_n(\{ n\} \times n),E)$ where $E$ is
equally of first coordinates). \newline
9) If $T_1$ is countable unstable, then $T = {\frac {\triangleleft^\ell}
{\ne}} T_{\text{eq}} {\frac {\triangleleft^\ell}{\ne}} T$ moreover
$\lambda > \aleph_0 \Rightarrow T_{\text{eq}} \triangleleft^*_\lambda T,
\lambda \ge 2^{\aleph_0} \Rightarrow T \triangleleft^*_\lambda T_{\text{eq}}$.
\enddemo
\bigskip

\demo{Proof}  1) - 4) Obvious. \newline
5) The proof of \cite{Sh:a},VI,2.6 = \cite{Sh:c},VI,2.6,p.337 gives this, too.
\newline
6), 7), 8)  As there [references]. \hfill$\square_{2.11}$
\enddemo
\bigskip

\demo{2.12 Proof of Theorem 2.9(1),(2)}  Without loss of generality
$\tau(T)$ is finite.  Remember: if $T'$ has infinite linear orders as models
then it is $\triangleleft^\ell$-maximal.  Let $J$ be a dense linear order,
such that:
\medskip
\roster
\item "{(a)}"  $J$ has a closed interval which is $I$
\item "{(b)}"  for any regular $\mu_1,\mu_2 \le |J|,J$ has an interval
isomorphic to $[(\{ 1\} \times \mu_1) \cup (\{2\} \times \mu_2)$ ordered by
$i_1,\alpha_1) <_J (i_2,\alpha_2) \Leftrightarrow (i_1 = 1 \and i_2 = 2)
\bigvee (i_1 = 1,i_2 \and \alpha_1 < \alpha_2) \bigvee (i_1 = 2 = \ell_2
\and \alpha_1 > \alpha_2)$.
\endroster
\medskip

\noindent
Let $\varphi(x,y)$ exemplify the SOP$_3$.  Let $M$ be a model of $T$ and
$F:I \rightarrow {}^{\ell g \bar x}M$ be such that $I \models \eta < \nu
\Rightarrow M \models \varphi [F(\eta),F(\nu)]$ and for every $c \in M^
{\frak B}$ we can find a finite $I' \subseteq I$ such that: if
$[t_1,t_2 \in (I \backslash I')] \and \dsize \bigwedge_{s \in I'}[s <_I t
\equiv s <_I t_2]$ then $M^{\frak B} \models \varphi[F(t_1),c] \equiv
\varphi[F(t_2),c]$ and $M^{\frak B} \models \varphi [c,F(t_1)] \equiv
\varphi [c,F(t_1)]$.  Let ${\frak B}_0 = (H(\chi),\in,<^*_\chi,J,I,F,M)$,
and ${\frak B}$ be a model, $\bold j$ an elementary embedding of ${\frak B}_0$
into ${\frak B}$ such that $M^* = M^{\frak B} \restriction L(T)$ is locally
$\lambda^+$-saturated but $I^{\frak B} = \bold j(I)$ is not $\lambda^+$-
saturated (for 2.9: ${\frak B}^* = {\frak B}^\lambda/D$).
\medskip

As $\bold j(I)$ is not $\lambda^+$-saturated, we can find $\lambda_0,
\lambda_1 \le \lambda$ and $\alpha^\ell_i \in \bold j(I)$ (for
$i < \lambda_\ell$, \newline
$\ell < 2$) such that:
\medskip
\roster
\item "{$(\alpha)$}"  $I^{\frak B} \models a^0_i < a^0_j$ for $i < j <
\lambda_0$
\item "{$(\beta)$}"  $I^{\frak B} \models a^1_i > a^1_i$ for $i < j <
\lambda_1$
\item "{$(\gamma)$}"  $I^{\frak B} \models a^0_i < a^1_j$ for $i < \lambda_0,
j < \lambda_1$
\item "{$(\delta)$}"  $I^{\frak B} \models \neg(\exists x)
[\dsize \bigwedge \Sb i < \lambda_0\\ j < \lambda_1 \endSb a^0_i < x <
a^1_j]$.
\endroster
\medskip

Clearly $\{ \varphi(\bar a^0_i,\bar x),\varphi(\bar x,\bar a^1_j): 
i < \lambda_0,j < \lambda_1\}$ is finitely satisfiable in $\bold j(M)$.  Now
as $\bold j(M)$ is locally $\lambda^+$-saturated there is $\bar a \in
{}^{\ell g \bar x}(M^*)$ such that $M \models$ ``$\varphi(\bar a^0_i,\bar a)
\and \varphi(\bar a,\bar a^1_j)$ for $i < \lambda_0, j < \lambda_1$".  In
${\frak B}$ we can define:

$$ 
I^{\frak B}_-[\bar a] =: \biggl\{ \eta \in I^{\frak B}:\text{there is }
\nu \in I^{\frak B} \text{ such that}: I^{\frak B} \models ``\eta \le \nu"
\text{ and } \bold j(M) \models \varphi(\bar a_\nu,\bar a) \biggr\}
$$

$$
\align
I^{\frak B}_+[\bar a] =: \biggl\{ \eta \in I^{\frak B}:&\text{there is }
\nu \in I^{\frak B} \text{ such that}: I^{\frak B} \models ``\nu \le \eta" \\
  &\text{and } j(M) \models \varphi[\bar a,\bar a_\nu] \biggr\}.
\endalign
$$
\medskip

\noindent
Clearly
\medskip
\roster
\item "{(a)}"  $I^{\frak B}_-[\bar a]$ is an initial segment of $\bold j(I)$
which belongs to ${\frak B}$.
\item "{(b)}"  $I^{\frak B}_+[t]$ is an end segment of $\bold j(I)$ which
belongs to ${\frak B}$.
\item "{(c)}"  For every $i < \lambda_0$
{\roster
\itemitem{ $(*)^i_0$ }  $a^0_i \in I^{\frak B}_-[\bold j,t]$
\endroster}
\item "{(d)}"  for every $j < \lambda_1$
{\roster
\itemitem{ $(*)^j_1$ }  $a^1_j \notin I^{\frak B}_-[t]$
\endroster}
\item "{(e)}"  By the choice of $\varphi$
{\roster
\itemitem{ $(*)_3$ }  $I^{\frak B}_-[t] \cap I^{\frak B}_+[t] = \emptyset$.
\endroster}
\endroster
\medskip

If for some $c \in {\frak B},{\frak B} \models ``c \in \bold j(J)$" and
$(\forall x \in I_-[t]) < t \le_I c)$ and \newline
$(\forall x \in I_+[t])[c \le_I x]$
we are done.  So ${\frak B}$ thinks $I^{\frak B}_-[t'],I^{\frak B}_+[t])$ is
a Dedekind cut, so let ${\frak B} \models ``\text{cf}(I^{\frak B}_-[t],<_I)
=t_-,\text{cf}(I^{\frak B}_+[t] >_I) = t_+$ and the (outside) cofinalities
of $t_-,t_+$ are $\mu_1,\mu_2$ respectively.  If $\mu_1,\mu_2$ are infinite,
we use clause (b) of the choice of $J$ (and the choice of $\mu$).  We are
left with the case where $\mu_1 = 1 < \mu_2$ (the other case is the same).
Use what ${\frak B}$ ``thinks" is a $(t_1,t_2)$ Dedekind cut of $J$ to show
$\mu_2 \ge \mu^+$ a contradiction. \hfill$\square_{2.9}$
\enddemo
\bigskip

\proclaim{2.13 Theorem}  1) The theorems on non-existence of a universal
model in $\lambda$ for linear order from \cite{KjSh:409}, 
\cite{Sh:457},\S3 hold
for any theory with SOP$_4$. \newline
2)  We can use embedding (not necessarily elementary) if $\varphi(\bar x,
\bar y)$ is quantifier free or even existential.
\endproclaim
\bigskip

\demo{Proof}  We concentrate on the case $\lambda$ is regular and part (1).
We will concentrate on the new part relative to \cite{KjSh:409}.  Let
$\varphi(\bar x,\bar y)$ exemplify SOP$_{\le 4}$ (exists by 2.6(1)) in a
complete first order theory $T$.  Without loss of generality
$\ell g \bar x = \ell g \bar y = 1$ and $T \vdash \neg \varphi
(\bar x,\bar x)$.
\medskip

Let $M$ be a model of $T$ with universe $\lambda,I$ a linear order,
$a_s \in M,M \models \varphi[a_s,a_t]$ for $s <_I t$ (from $I$).

We do not have a real Dedekind cut (as $\varphi(x,y)$ is not transitive),
but we use replacements.  Now for every $b \in M$, let $I^-[b] = \{ t:M
\models \varphi [a_t,b]\}$ and \newline
$I^+[b] = \{ t:M \models \varphi[b,a_t]\}$.
As $\varphi$ exemplifies also SOP$_{\le 3}$ clearly the following is
satisfied:
\medskip
\roster
\item "{$(*)$}"  $s \in I^-[b] \and t \in I^+[t] \Rightarrow s < t$ \newline
(if $t < s$ from a counterexample, $b,t,s$ gives a contradiction).
\endroster
\medskip

\noindent
Note: $I^-[a_s] = \{ t:t <_I s\}, I^+[a_s] = \{ t:s <_I t\}$. \newline
Let $P = \{ a_s:s \in I\},<^* = \{(a_s,a_t):s <_I t\}$ and \newline
$J^-[t,a] = \{ s \in I:s <_I t,a_s \in \alpha\},J^+[t,\alpha] =
\{ s \in I:t <_I s,a_s \in \alpha\}$ \newline
(remember: $|M|$, the universe of $M$,
is $\lambda,\alpha = \{ \beta:\beta < \alpha\}$).  \newline
Hence
$C =: \{ \delta < \lambda:(M,P,<^*,<) \restriction \delta \prec M^+ =:
(M,P,<^*,<)\}$ is a club of $\lambda$.  Clearly:
\medskip
\roster
\item "{$(**)$}"  let $\delta \in C,b \in M \cap \delta$;
{\roster
\itemitem{ (i) }  if $(I^-[b],<_I)$ has cofinality $< \lambda$, 
\underbar{then} $I^-[b] \cap M_\delta$ is $<_I$-cofinal in it
\itemitem{ (ii) }  if $(I^+[b],>_I)$ has cofinality $< \lambda$, 
\underbar{then} $I^+[b] \cap M_\delta$ is $(>_I)$-cofinal in it
\itemitem{ (iii) }  if there is $t$ such that $I^-[b] \le_I t \le_I I^+[b]$
\underbar{then} there is such $t \in M \cap \delta$.
\endroster}
\endroster
\medskip

\noindent
Now suppose that $\delta_1 < \delta_2$ are in $C,t(*) \in I,a_{t(*)} \in P
\backslash \delta_2$.
\enddemo
\bigskip

\noindent
\underbar{Case 1}:  For some $s(*) \in M \cap \delta_2$ we have
$(\forall s \in J^-[t(*),\delta_1])(s <_I s(*) <_I t(*))$.  Let $b =:
a_{s(*)}$.  Hence for every $ c \in M \cap \delta_1$: if $\varphi(c,a_{t(*)})$
then for every $t',t''$ satisfying $t" <_I t(*) <_I t'',t' \in I,a_{t'} \in
\delta_1,t'' \in I,a_{t''} \in \delta_1$ we have $M^+ \models (\exists x)
[x \in P \and \varphi(c,x) \and a_{t'} <^* x <^* <a_{t''}]$.  Clearly (or see
the middle of the proof of case 2 below) necessarily for arbitrary $<_I$-
large $t \in J^-[t(*),\delta_1]$ we have $\varphi[c,a_t]$ but for any such
$t,\varphi[a_t,a_{s(*)}]$ i.e. $\varphi[a_t,b]$ hence
\medskip
\roster
\item "{$(*)_1$}"   $(\forall c \in M \cap \delta_1)[\varphi(c,a_{t(*)})
\rightarrow (\exists y \in M \cap \delta_1)[\varphi(c,y) \and \varphi(y,b)]$.
\endroster
\medskip

\noindent
Of course,
\medskip
\roster
\item "{$(*)_2$}"  $b \in M \cap \delta_2$
\item "{$(*)_3$}"  $\varphi[b,a_{t(*)}]$.
\endroster
\medskip

\noindent
Note: those three properties speak on $M,\delta_1,\delta_2,a_{t(*)},b$ but
not on $I,<^*,P,<$.
\bigskip

\noindent
\underbar{Case 2}:  For no $s(*) \in I,a_{s(*)} \in \delta_2$ do we have
$(\forall s \in J^-[t(*),\delta_1])[s <_I s(*) <_I t(*)]$ we assume:
\medskip
\roster
\item "{(A)}"  $\{ a_s:s \in I,a_s \in \delta_2,s <_I t(*)\}$ is not definable
in $M^+ \restriction \delta_2$.
\endroster
\medskip

We shall now show that for no $b \in M \cap \delta_2$ do we have $(*)_1 +
(*)_2 + (*)_3$, so assume $b$ is like that and we shall get a contradiction.

Without loss of generality $(A)$ holds.  By $(*)_3$ we have $\varphi[b,
a_{t(*)}]$ hence for arbitrarily $<_I$-large $t \in J^-[t(*),\delta_2]$ we
have $\varphi[b,a_t]$; choose such $t_0$. \newline
[Why?  Otherwise $I^+[b] \cap J^-[t(*),\delta_2]$ is bounded say by some
$t^*$, so $\theta(x,b,a_{t^*}) =: x \in P \and (\exists y)[y \in P \and y
\le^* x \and \varphi(b,y) \and t^* <^* x]$ define in $M^+$ a set which is an
end segment of $(P,<^*)$, include $t(*)$ (check) but so $s \in \delta_2,
s <^* t(*)$.  So in $M^+ \restriction \delta_2$ it defines the set
$\{ a_s:s \in J^+[t(*),\delta_2]\}$ hence $\neg \theta(x,b,a_{t^*})$ define
in $M^+ \restriction \delta_2$ the set $\{ a_s:s \in J^-[t(*),\delta_2]\}$,
hence by the assumption of the case, $M^+ \restriction \delta_2$ satisfies:

$$
(\forall z)[z \in P \and \neg \theta(z,b,a_{t^*}) \rightarrow (\exists y <
\delta_1)(y \in P \and z \le^* y \and \neg \theta(z,b,a_{t^*})]
$$
\medskip

\noindent
contradicting $(A)$ above. \newline

So by the assumption of the case (i.e. that $t_0 < \delta_1$ cannot serve
as $s(*)$ and $t_0 <_I t(*)$) for some $t_1 \in J^-[t(*),\delta_1]$ we have
$t_0 <_I t_1$ and clearly $t_1 <_I t(*)$ hence $\varphi[a_{t_1},a_{t(*)}]$.
So by $(*)_1$ applied with $a_{t_1}$ standing for $c$ for some $y \in M \cap
\delta_1$ we have $\varphi[a_{t_1},y] \and \varphi[y,b]$.  Now $b,a_{t_0},
a_{t_1},y$ contradicts ``$\varphi(x,y)$ exemplifies SOP$_4$". \newline
Hence we get together
\medskip
\roster
\item "{$\bigoplus$}"  if $\delta_1 < \delta_2 \in C,a_{t(*)} \in P \backslash
\delta_2$ (so $t(*) \in I)$ then the following conditions are equivalent:
{\roster
\itemitem{ $(\alpha)$ }  for some $s(*) \in \delta_2$, (so necessarily
$s(*) \ne t(*)$) we have $(\forall s \in I)[a_s \in \delta_1 \Rightarrow s
<_I s(*) \equiv s <_I t(*)]$) we have: \newline

$\qquad (A) \quad \{ a_s:s \in I,a_s \in \delta_2 s <_I t(*)\}$ 
is definable in $M^+ \restriction \delta_2$ \newline

$\qquad \qquad$ (with parameters)
\newline

$\qquad (B) \quad \{ a_s:s \in I,a_s \in \delta_2,t(*) <_I s\}$ is
definable in $M^+ \restriction \delta_2$ \newline

$\qquad \qquad$ (with parameters)
\smallskip
\noindent
\itemitem{ $(\beta)$ }  for some $b \in \delta_2$ the conditions $(*)_1,
(*)_2,(*)_3$ above holds for $\varphi(x,y)$ or for $\varphi^-(x,y)$ where
$\varphi^-(x,y) = \varphi(y,x)$.
\endroster}
\endroster
\bigskip

\demo{Proof}  If clause $(\alpha)$ holds and $s(*) <_I t(*)$ holds use Case
1 above.  If clause $(\alpha)$ holds and $s(*) <_I t(*)$ fails, then
$t(*) <_I s(*)$ inverts the order of $I$, use $\varphi^-$ and now apply Case
1 above.  So assume $\neg(\alpha)$.  We first want to apply Case 2 to prove
there is no $b$ satisfying $(*)_1,(*)_2,(*)_3$.  For this we need clause (A)
there.  We claim it holds. \newline
[Why?  Assume $\bar d \in (M^2 \restriction \delta_2),\psi$ a first order
formula (in the vocabulary of $M^+$), such that for every $e \in M^+ 
\restriction \delta_2$ we have: $M^+ \restriction \delta_2 \models \psi[e,
\bar d]$ iff $e \in \{ a_s:s \in I,a_s \in \delta_2,s <_I t(*)\}$.  So
$M^+ \models (\exists z)[P(z) \and (\forall y)(y < \delta_1 \and P(y) 
\Rightarrow y <^* z \equiv \psi[y,\bar d])]$ as $z \mapsto a_{t(*)}$ satisfies
it, but $M^+ \restriction \delta_2 \prec M^+$ hence there is $z^* \in 
\delta_2$ satisfying this.  So $z^* \in P$ hence for some $s(*)$, $z^* =
a_{s(*)}$, so $s(*)$ contradicts the assumptions $\neg(\alpha)$.  So we have
proved the failure of the first possibility from clause $(\beta)$.  The
second is proved similarly inverting the order of $I$, using $\varphi^-$
(noting that this transformation preserves the statement (A) from Case 2].
\enddemo
\bigskip

\definition{2.14 Definition}  Let $M$ be a model with universe $\lambda$ and
$\varphi(x,y)$ a formula exemplifying SOP$_4$ (possibly with parameters)
let $\varphi^+(x,y) =: \varphi(x,y),\varphi^-(x,y) =: \varphi(y,x)$.  Assume
$\bar C = \langle C_\delta:\delta \in S \rangle$ is a club system, $S
\subseteq \lambda$ stationary, guessing club \footnote {otherwise dull} (i.e.
for every club $E$ of $\lambda$ for stationarily many $\delta < \lambda,\delta
\in S,C_\delta \subseteq E$)
\medskip
\roster
\item "{(a)}" for $x \in |M|$ and $\delta \in S$ let
$$
\align
\text{inv}_\varphi(x,C_\delta,M) = \biggl\{ 
\alpha \in \text{ nacc }C_\delta:&\text{letting }\delta_2 = \alpha,
\delta_1 = \sup(C_\delta \cap \alpha) \text{ (well defined)}, \\
  &\text{for some } b \text{ conditions }
(*)_1,(*)_2,(*)_3 \text{ of Case 1} \\
  &\text{holds for } \varphi^+ \text{ or for } \varphi^- \biggr\}
\endalign
$$
\item "{(b)}"  $\text{Inv}_\varphi(C_\delta,M) = \{\text{inv}(x,C_\delta,M):
x \in M\}$ \newline
$\text{INv}_\varphi(M,\bar C) = \langle \text{Inv}_\varphi(C_\delta,M):\delta
\in S \rangle$ \newline
$\text{INV}_\varphi(M<\bar C) = \text{ INv}_\varphi(M,\bar C)/\text{id}^\alpha
(\bar C)$ where:
\endroster
\enddefinition   
\bigskip

\definition {2.15 Definition}  \newline
$\text{id}^\alpha(\bar C) = \{S' \subseteq \lambda:\text{for some club } E
\text{ of } \lambda \text{ the set of } \delta \in S' \cap S$ for which
$C_\delta \subseteq E$ is not stationary$\}$.
\enddefinition
\bigskip

\demo{2.16 Observation}  If $M' \cong M''$ are models of $T$ and both have
universe $\lambda$ in $M$ then $\text{INV}_\varphi(M',\bar C) = \text{ INV}
_\varphi(M'',\bar C)$ so $\text{INV}_\varphi(M,\bar C)$ can be defined for
any model of cardinality $\lambda$.
\enddemo
\bigskip

\demo{Proof}  Let $f$ be from $M'$ onto $M''$, so $f$ is a permutation of
$\lambda$.  So $E_0 = \{ \delta < \lambda:\delta \text{ a limit ordinal}, f
\text{ maps } \delta \text{ onto } \delta\}$.  Assume $C_\delta \subseteq
E$, then for $x \in M' \backslash \delta,\delta \in S$ we have inv$(x,
C_\delta,M') = \text{ inv}_\varphi(f(x),C_\delta,M'')$. \newline
[Why?  Read $(*)_1,(*)_2,(*)_3$].  Hence Inv$_\varphi(C_\delta,M') \in 
\text{ Inv}_\varphi(C_\delta,M'')$.  By the definition of $\text{id}^\alpha
(\bar C)$ we are done. \hfill$\square_{2.16}$
\enddemo
\bigskip

\demo{2.17 Observation}  If $\bar C = \langle C_\delta:\delta \in S \rangle,
S \subseteq \lambda$ stationary, $C_\delta \subseteq \delta = \sup(C_\delta),
C_\delta$ closed, $I$ a linear order with the set of elements being $\lambda$
we let:
\medskip
\roster
\item "{(a)}"  for $x \in \lambda,\delta \in S,\text{inv}(x,C_\delta,I) =
\{ \alpha \in \text{ nacc}(C_\delta):\text{there are } y,z \in \alpha
\text{ such}$ \newline
$\text{that } y <_I x <_I z \text{ such that } (\forall s)[s < \sup
(C_\delta \cap \alpha) \Rightarrow s <_I y \vee z <_I s]\}$
\item "{(b)}"  Inv$(C_\delta,I) = \{\text{inv}(x,C_\delta,I):x \in M\}$
\item "{(c)}"  INv$(I,\bar C) = \langle\text{Inv}(C_\delta,I):\delta \in S
\rangle$
\item "{(d)}"  INv$(I,\bar C) = \text{ INv}(I,\bar c)/\text{id}^\alpha
(\bar C)$
\endroster
\enddemo
\bigskip

\demo{2.18 Observation}  INV$(I,\bar C) = \text{ INV}(I',\bar C)$ if $I \cong
I'$, so actually it is well defined for any linear order with cardinality
$\lambda$.
\enddemo
\bigskip

\demo{2.19 Observation}  If $M$ is a model with universe $\lambda$ and
$\varphi,\langle a_s:s \in I \rangle$ as above and, $\emptyset \notin
\text{ id}^\alpha(\bar C)$ \underbar{then} INV$(I,\bar C) \le \text{ INV}
_\varphi(M,\bar C)$ i.e. for some club $E$ of $\lambda,\delta \in S \and
C_\delta \subseteq E \Rightarrow \text{ Inv}(C_\delta,I) \subseteq
\text{ Inv}_\varphi(C_\delta,M)$.
\enddemo
\bigskip

\demo{Proof}  By $\oplus$ above.
\enddemo
\bigskip

\demo{Conclusion of the proof of 2.13}  As in \cite{KjSh:409}. \hfill
$\square_{2.13}$
\enddemo
\bigskip

\proclaim{2.20 Claim}  For a complete $T$, the following are equivalent:
\medskip
\roster
\item "{(a)}"  $T$ does not have SOP$_3$
\item "{(b)}"  if in ${\frak C},\langle \bar a_i:i < \alpha \rangle$ is an
indiscernible sequence, $\alpha$ infinite and \newline
$\{ \varphi(\bar x,\bar y,\psi
(\bar x,\bar y)\}$ contradictory and for each $j$ for some $\bar b_j$ we
have \newline
$i \le j \Rightarrow \models \varphi[\bar a,\bar b,\bar a_1]$ and
$i > j \Rightarrow \models \psi[\bar b,\bar a_i]$ \underbar{then} for 
$i < j$ \newline
we have 
$(\exists \bar x)(\varphi(\bar x,\bar a_j) \and \psi(\bar x,\bar a_i))$
\item "{(c)}"  in clause (b) we replace the conclusion: for every finite
disjoint $u,v \subseteq \omega$ \newline
we have
$(\exists \bar x)\left( \dsize \bigwedge_{i \in u} \varphi(\bar x,\bar a_i)
\and \dsize \bigwedge_{j \in v} \psi(\bar x,\bar a_j) \right)$.
\endroster
\endproclaim
\bigskip

\demo{Proof}  \underbar{$(c) \Rightarrow (b)$}:  Trivial. \newline
\medskip
\noindent
\underbar{$\neg(c) \Rightarrow \neg(b)$}:  Choose counterexample 
with $|u \cup v|$  minimal, assume $\alpha > \omega + |u \cup v|$. \newline
\medskip
\noindent
\underbar{$\neg(a) \Rightarrow \neg(b)$}:  
Straight by the Definition of SOP$_3$, etc.
\newline
\medskip
\noindent
\underbar{$\neg(b) \Rightarrow \neg(a)$}:  Without loss of generality
$\langle \bar a_i \char 94 \bar b_i:i < \alpha \rangle$ is an indiscernible
sequence.  Now we cannot find $\bar c_0,\bar c_1,\bar c_2$ such that
$\bar c_0 \char 94 \bar c_1, \bar c_1 \char 94 \bar c_2,\bar c_2 \char 94
\bar c_0$ realizes the same type as $(\bar a_0 \char 94 \bar b_0) \char 94
(\bar a_1 \char 94 \bar b_1)$, so SOP$_3$ is exemplified. \hfill$\square_{
2.20}$
\enddemo
\newpage

\nocite{HrPi1}
\nocite{HrPi2}
\nocite{ChHr}

\bigskip
REFERENCES
\bigskip

\bibliographystyle{literal-plain}
\bibliography{lista,listb,listx,listf}
\shlhetal

\enddocument
\bye